%% file: framesAMcorrected2.tex
\documentclass[11pt, reqno]{amsart}
\input{macros}

\usepackage{amsmath,amsthm,amsfonts,amssymb}
\usepackage{epsf}
\usepackage{epsfig}

\newcommand{\simca}{\sim _{\cA }}

\title{Banach frames for $\alpha$-modulation spaces}
\author{Massimo\ Fornasier }
\begin{document}
\maketitle
\begin{center}
Dipartimento di Metodi e Modelli Matematici\\
per le Scienze Applicate\\
Universit\`a di Roma ``La Sapienza''\\
Via A. Scarpa 16/B\\
I-00161 Roma, Italia\\
and\\

NuHAG, Department of Mathematics\\
University of Vienna\\
Nordbergstrasse 15\\
A-1090 Vienna, Austria
\bigskip

email: {\tt mfornasi@math.unipd.it}
\end{center}

\begin{abstract}
This paper is concerned with the characterization of $\alpha$-modulation spaces by Banach frames, i.e., stable and redundant non-orthogonal expansions, constituted of functions obtained by a suitable combination of translation, modulation and dilation of a mother atom. In particular, the parameter $\alpha \in [0,1]$ governs the dependence of the dilation factor on the frequency. The result is achieved by exploiting intrinsic properties of localization of such frames.
The well-known Gabor and wavelet frames arise as special cases  ($\alpha = 0$) and limiting case ($ \alpha \to 1)$, to characterize respectively modulation and Besov spaces. This intermediate theory contributes to a further answer to the theoretical need of a common interpretation and framework between Gabor and wavelet theory and to the construction of new tools for applications in time-frequency analysis, signal processing, and numerical analysis. 
\end{abstract}



\section{Introduction}
The theory of {\it frames}, or stable redundant non-orthogonal
expansions in Hilbert spaces, introduced  by Duffin and Schaeffer \cite{DS},
plays an important role in  {\it wavelet theory} \cite{D1,D2,DGM} as well as
in {\it Gabor (time-frequency) analysis} \cite{Gr,FS,FS1} for functions in $L^2(\mathbb{R}^d)$.
Besides traditional and relevant applications of frames in signal processing, image processing, data compression,   pattern matching, sampling theory, communication and data transmission, recently the use of frames also in numerical analysis for the solution of operator equation is investigated \cite{St,DFR}.
Therefore, not only the characterization  by frames of functions in $L^2(\mathbb{R}^d)$ is relevant but also that of (smoothness) Banach  function spaces is crucial to have a correct formulation of effective and stable numerical schemes. The concept of {\it Banach frame} as an extension of atomic decompositions in {\it coorbit spaces} \cite{FG1,FG2} has been already introduced in \cite{gro}.
Moreover this classical theory of Feichtinger and Gr\"ochenig has shown in particular that  Gabor and wavelet $L^2$-frames can in fact extend to Banach frames for {\it modulation} \cite{F4,Gr,gro1,FoG} and {\it (homogeneous) Besov spaces} \cite{FJ,Tr,Tr2} respectively. 
As a further answer to the theoretical need of a common interpretation and framework between Gabor and wavelet theory, the author has recently proposed \cite{FF} the construction of frames,
which allows to ensure that certain families of Schwartz functions (atoms) on $\mathbb{R}$ obtained
by a suitable combination of translation, modulation and dilation 
\begin{eqnarray*}
       && T_x ( f )(t) = f(t-x),\\ 
       && M_\omega ( f )(t) = e^{2 \pi i \omega \cdot t} f(t),\\
       & & D_a ( f )(t) = a^{-1/2} f(t/a), \quad x, \omega,t \in \mathbb{R}, a \in \mathbb{R}_+,\\
\end{eqnarray*}
form Banach frames for the family of $L^2$-Sobolev spaces of any order.
In this construction a parameter $\alpha \in [0,1)$ governs the
dependence of the dilation factor 
on the frequency parameter.
The well-known Gabor and wavelet frames  (also valid for the same
scale of Hilbert spaces that constitutes an intersection of the modulation and Besov space families) 
arise as special cases  ($\alpha = 0$) and limiting case ($ \alpha \to 1)$ respectively.  Thus, let us call these families $\alpha$-Gabor-wavelet frames. In contrast to those
limiting cases it is no longer possible to use group theoretic
arguments nor the coorbit space theory can be applied anymore to extend the $L^2$-frame to a Banach frame. A similar approach was proposed by Hogan and Lakey \cite{HL} to construct {\it coherent frames} generated by representations of extensions of the Heisenberg group by dilation. Other contributions due to Weiss {\it et al.}  \cite{HLW,HLWW,L} developed characterizations of a large class of mixed decompositions in $L^2$ as an attempt of a unified approach to Gabor, wavelet, and more general wave packet frames. 

New tools for extending an $L^2$-frame to Banach frames and atomic decompositions have been introduced by Gr\"ochenig.
The key concept in \cite{gro1} is the localization properties of the frame with respect to an auxiliary Riesz basis. The localization has been measured by polynomial or sub-exponential off-diagonal decay of the cross Gramian matrix of the frame and the Riesz basis. The main result in \cite{gro1}  asserts that a localized frame has canonical dual with the same localization properties and that the frame extends to a Banach frame and an atomic decomposition for the Banach spaces for which the reference auxiliary Riesz basis is a unconditional basis. 
Inspired by this work, the author \cite[Chapter 5]{Fo2} showed that the extension of a frame to Banach frames does not depend on localization properties with respect to any auxiliary Riesz basis, but it can be formulated also as an intrinsic property of the frame. In particular, if the frame is {\it intrinsically or self-localized}, i.e., if its Gramian matrix has a suitable off-diagonal decay, and there exists a corresponding {\it dual frame} with the same property then the frame extends in fact to a Banach frame and an atomic decomposition for a suitable class of Banach spaces.
Based on a rather tricky and technical construction of an intrinsically localized dual frame,  this principle has been applied in \cite[Chapter 5]{Fo2} to extend $\alpha$-Gabor-wavelet $L^2$-frames to atomic decompositions for {\it $\alpha$-modulation spaces}. These Banach (smoothness) function spaces have been introduced independently by Gr\"obner \cite{G} and Paiv\"arinta/Somersalo \cite{PS} as an ``intermediate'' family between modulation and Besov spaces. They appear also as particular cases of the spaces introduced by Holschneider and Nazaret in \cite [Section 4.2]{HN}, and Hogan and Lakey in \cite[Section 4.5]{HL2}, by retract or pull back methods based on generalized {\it Fourier-Bros-Iagolnitzer transforms} \cite{BI} (or {\it flexible Gabor-wavelet transforms} as they are called in \cite{FF,Fo2}). Characterizations of $\alpha$-modulation spaces by {\it brushlet unconditional basis} have been given by Nielsen and Borup \cite{NB} and the mapping properties of pseudodifferential operators in H\"ormander classes on $\alpha$-modulation spaces have been studied by  Holschneider and Nazaret \cite{HN} and   Borup \cite{B}, as generalizations of classical results of Cordoba and Fefferman \cite{CF}.

In this paper we shall present a Banach frame and atomic decomposition characterization of $\alpha$-modulation spaces, following the intrinsic localization strategy already suggested in \cite[Chapter 5]{Fo2}.
The result will be achieved firstly by describing functions in  $\alpha$-modulation spaces by means of suitable families of band-limited functions, and then extending the result to $\alpha$-Gabor-wavelet frames by means of general perturbation principles, here applied exploiting localization properties of such frames.

The paper is organized as follows. Section 2 recalls the concept of frames in Hilbert and Banach spaces. In particular, the intrinsic localization of frame theory is discussed as a method to extend frames in Hilbert spaces to Banach frames.
In Section 3 we present  $\alpha$-modulation spaces as a generalization of modulation and inhomogeneous Besov spaces and the localization principles applied to $\alpha$-Gabor-wavelet frames to characterize them. We conclude with few remarks and a characterization of $\alpha$-modulation spaces by pull back of certain weighted $L^{p,q}$ spaces (mixed norm Lebesgue spaces) by the {\it flexible Gabor-wavelet transform} introduced in \cite{HN,FF,Fo2}.
\\

 \emph{Acknowledgment:} 
The author thanks Hans G. Feichtinger and Karlheinz Gr\"ochenig for the fruitful discussions, their valuable suggestions and the hospitality of NuHAG (the Numerical Harmonic Analysis Group, Department of Mathematics, University of Vienna, AUSTRIA) during the preparation of this work.

A special thank is addressed to both the anonymous referees and Joachim St\"ockler for the careful reading of the manuscript and for the relevant suggestions in order to improve its quality.
 The author acknowledges the support of the Intra-European Individual Marie Curie Fellowship, project FTFDORF-FP6-501018.

\subsection{Notations}
We denote with $L^p(\mathbb{R}^d)$ the Lebesgue space of measurable functions on $\mathbb{R}^d$ that are $p$-integrable and with $L^p_m(\mathbb{R}^d)$  the Lebesgue space of measurable functions $f$ such that $f m \in L^p (\mathbb{R}^d)$. Similarly are defined the spaces $\ell^p_m(\mathbb{Z}^d)$ of weighted $p$-summable sequences. The space $\mathcal{S}(\mathbb{R}^d)$ is the space of Schwartz functions and its dual $\mathcal{S}'(\mathbb{R}^d)$ is the space of tempered distributions. We denote with $\mathcal{F}$ the Fourier transform on $\mathcal{S}'(\mathbb{R}^d)$ and with $\mathcal{F} L^p$ the space of distributions which are images of $L^p$ functions under the action of $\mathcal {F}$, endowed with the natural norm $\|f\|_{\mathcal{F} L^p}:=\|\mathcal F^{-1} f\|_p$. For positive quantities $F$ and $G$, we will write $F \lesssim G$ whenever $F(x) \leq C \cdot G(x)$ for some universal constant $C>0$ and for all variable $x$. When $F \lesssim G$ and $G \lesssim F$ then we will write $F \asymp G$. For any function $g$ on $\mathbb R$ we define the operator $\cdot^\nabla$ by $g^\nabla(t):=g(-t)$, for all $t \in \mathbb R$. The function $\text{sgn}(x)=1$ if $x>0$, $\text{sgn}(x)=-1$ if $x<0$, and $\text{sgn}(x)=0$ if $x=0$. The symbol $\chi_E$ denotes the characteristic function of $E \subset \mathbb{R}$.

\section{Intrinsically localized frames in Banach spaces}

\subsection{Frames in Hilbert and Banach spaces} 
In this section we recall the concept of frames, how they  can be used to define certain associated
Banach spaces, and how to obtain  stable decompositions in these  Banach spaces.


A subset  $\cG=\{g_n\}_{n \in {\mathbb{Z}}^d}$ of a separable Hilbert
 space   $\mathcal{H}$ is called \emph{frame} for $\mathcal{H}$ if
\begin{equation}
 A \|f\|^2 \leq \sum_{n \in {\mathbb{Z}}^d} |\langle f, g_n \rangle|^2 \leq B \|f\|^2, \quad \forall f\in \mathcal{H},
\end{equation}
for some  constants $0<A \leq B < \infty$.

Equivalently, we could define a frame by the requirement that  the
corresponding  \emph{analysis   operator} $C=C_{\cG }$ defined by
$Cf = (\langle f, g_n \rangle)_{n \in \mathbb{Z}^d} $ is bounded from $\mathcal{H}$ into $\ell^2({\mathbb{Z}}^d)$ or  that the  \emph{synthesis operator} $D=
D_{\cG } = C^*, D\mathbf{c} = \sum_{n \in \mathbb{Z}^d} c_n g_n, $  is bounded from $
\ell^2({\mathbb{Z}}^d)$ into $\mathcal{H}$, and the frame
operator  $S=D  C$ is  boundedly invertible (positive and
self-adjoint) on $\cH $. The family  $\tilde{\cG}=S^{-1} \cG:=\{S^{-1} g_n\}_{n \in \mathbb{Z}^d}$ is again a
frame for $\mathcal{H}$. This  so-called  \emph{canonical dual frame}
plays an important role in the reconstruction of $f\in \cH $ from the
frame coefficients and in non-orthogonal expansions, because we have
\begin{equation}
f =S S^{-1} f = \sum_{n \in \mathbb{Z}^d} \langle f, S^{-1} g_n \rangle g_n=S^{-1} S f = \sum_{n \in \mathbb{Z}^d} \langle f, g_n \rangle S^{-1} g_n.
\end{equation}
Since in general  a frame is  overcomplete,  the coefficients in this
expansion
are  in general not unique (unless $\cG $ is a Riesz basis, we have
$\text{ker}(D) \neq \{0\}$) and  there may exist many possible other dual frames  $\{
\tilde g_n \}_{n \in {\mathbb{Z}}^d}$ in $\mathcal{H}$ such that
$$
        f = \sum_{n \in {\mathbb{Z}}^d} \langle f, \tilde g_n \rangle g_n
$$
with the norm equivalence $\|f\|_{\cH } \asymp \|\langle f, \tilde g_n \rangle_{n \in \mathbb{Z}^d}
\|_{\ell^2}$.
More information on frames can be found in the book \cite{chr}.
The concept of frame can be extended to Banach spaces as follows:

\begin{definition}
A \emph{Banach frame} for a separable Banach space $B$ is a sequence $\cG =\{g_n\}_{n \in {\mathbb{Z}}^d}$ in $B'$ with an associated sequence space $B_d$ such that the following properties hold.
\begin{itemize}
\item[(a)] The \emph{coefficient operator} $C$ defined by $C f=\L \langle f,g_n\rangle_{n \in {\mathbb{Z}}^d} \R$ is bounded from $B$ into $B_d$.
\item[(b)] Norm equivalence:
$$
        \|f\|_B \asymp \| \langle f,g_n\rangle_{n \in {\mathbb{Z}}^d}\|_{B_d}.
$$
\item[(c)] There exists a bounded operator $R$ from $B_d$ onto $B$, a
  so-called \emph{synthesis or reconstruction operator}, such that
$$
        R\L \langle f,g_n\rangle_{n \in {\mathbb{Z}}^d}\R =f.
$$
\end{itemize}
\end{definition}

As a dual concept and a different extension of Hilbert frames to Banach spaces is the notion of {\it atomic decomposition}.
\begin{definition}
 A \emph{frame atomic decomposition} for a separable Banach space $B$ is a sequence $\cG =\{g_n\}_{n \in \mathbb{Z}^d}$ in $B$ with an associated sequence space $B_d$ such that the following properties hold.
\begin{itemize}
\item[(a)] There exists a \emph{coefficient operator} $C$ defined by $C f=\L \langle f,\tilde g_n\rangle_{n \in\mathbb{Z}^d} \R$ bounded from $B$ into $B_d$, where $\tilde \cG =\{\tilde g_n\}_{n \in\mathbb{Z}^d}$ is in $B'$;
\item[(b)] norm equivalence:
$$
        \|f\|_{B} \asymp \| \langle f,\tilde g_n\rangle_{n \in\mathbb{Z}^d}\|_{B_d};
$$
\item[(c)] the following series expansion converge unconditionally
$$
        f= \sum_{n \in \mathbb{Z}^d} \langle f, \tilde g_n\rangle g_n, \quad \text{for all } f\in B.
$$
\end{itemize}
\end{definition}

In the following we discuss under which (sufficient) conditions and  for which suitable associated Banach
spaces  a Hilbert frame is also a Banach frame and an atomic decomposition.
In particular, this problem has motivated the theory of localized frames recently introduced by Gr\"ochenig \cite{gro1,GL,cg,FoG}.

\subsection{Intrinsic localization of frames} \label{localization}

We want to recall here the concept of mutual localization of two frames measured by their (cross-)Gramian matrix belonging to a 
class $\cA$ of matrices 
with suitable off-diagonal decay and mapping properties. The theory of localized frames 
has been introduced in \cite{gro1,GL} and recently developed in \cite{cg,FoG,FR}. In particular in case $\cA$ is a {\it spectral Banach $*$-algebra}   it has been shown  that a localized frame can extend to a Banach frame in a natural way for a large family of Banach spaces together with its canonical dual. We refer to  \cite{GL,FoG} for further information where a characterization of a large class of algebras of this type is presented.

In this paper we shall work with classes of matrices which are not necessarily algebras. As we will see, this will arise significant technical difficulties for the characterization of Banach spaces, which we can solve only by the use of the auxiliary construction of simpler frames and the applications of suitable perturbation results \cite{CH}. 
In the following we require that 
\begin{itemize}
\item[(A0)] $\cA \subseteq \cB (\ell ^2({\mathbb{Z}}^d ) )$, i.e., each $A \in
  \cA $ defines a bounded operator on $\ell ^2 ({\mathbb{Z}}^d )$.
\item[(A1)] $\cA $ is solid:  i.e., if $A\in \cA $ and  $|b_{kl}| \leq
  |a_{kl}|$ for all $k,l\in {\mathbb{Z}}^d $,  then  $B\in \cA$ as well.
\end{itemize}
 
Let us denote $w_s(x) =(1+|x|)^s$, for $s\geq0$, the polynomially growing submultiplicative and radial symmetric weight function on $\rd$. A weight $m$ on $\rd$ is called $s$-moderate if $m(x+y) \leq w_s(x) m(y)$. In particular, if $m$ is $s$-moderate then $m^{-1}$ is also $s$-moderate and $m(x) \lesssim w_s(x)$ for all $x \in \rd$.
As an additional requirement for Banach spaces characterization, we also ask that any $A \in \cA$ extends to a bounded operator from $\ell^p_m$ to $\ell^p_m$, for $1 \leq p \leq \infty$ and for suitable $s$-moderate weights $m$. By means of the class $\cA$, we can now state the  general localization concept.

Given two frames $\cG=\{g_n\}_{n \in {\mathbb{Z}}^d}$ and $\mathcal{F}=\{f_x\}_{x \in {\mathbb{Z}}^d}$ for the Hilbert space $\mathcal{H}$, the (cross-) Gramian matrix $A=A(\cG,\cF)$ of $\cG $ with respect to $\cF $ is the $\zd \times {\mathbb{Z}}^d $-matrix with entries
$$
a_{nx} = \langle g_n, f_x \rangle .
$$
A frame $\cG $ for $\cH $ is called \textit{$\cA$-localized} with respect to
another frame $\cF $ if $A(\cG,\cF) \in \cA$. In this case we write $\cG \sim
_{\cA } \cF $.  If $\cG \simca \cG$, then $\cG$ is called \textit{$\cA $-self-localized} or
\textit{intrinsically $\cA$-localized}. 


%

\subsection{Associated Banach Spaces} \label{assobanach}
In this subsection, we want to illustrate how $\cA$-self-localized frames can characterize suitable families of Banach spaces in a natural way.
In the following we assume $s \geq 0$ and $m$ is an $s$-moderate weight, and that $\cA \ell^p_m \subset \ell^p_m$ continuously, for all $p\in [1,\infty]$.

Let $(\cG,\tilde \cG)$ be a pair of dual  $\cA$-self-localized frames for $\mathcal{H}$ with $\cG \sim_{\cA} \tilde \cG$.
Assume $\ell^p_m({\mathbb{Z}}^d) \subset \ell^2({\mathbb{Z}}^d)$. Then the Banach space $\mathcal{H}_m^p(\cG,\tilde \cG)$ is defined to be
\begin{equation}
\label{Hpwdef}
\mathcal{H}_m^p(\cG,\tilde{\cG}) := \{ f \in \mathcal{H} : \quad f= \sum_{n \in {\mathbb{Z}}^d} \langle f, \tilde g_n \rangle g_n, \quad (\langle f,\tilde g_n\rangle)_{n \in {\mathbb{Z}}^d} \in \ell^p_m({\mathbb{Z}}^d)\}
\end{equation}
with the norm $\|f\|_{\mathcal{H}_m^p}= \| (\langle f,\tilde g_n \rangle )_{n \in {\mathbb{Z}}^d}\|_{\ell^p_m}$ and $1 \leq p \leq \infty$. 
Since $\ell^p_m({\mathbb{Z}}^d) \subset \ell^2({\mathbb{Z}}^d)$, $\mathcal{H}_m^p$ is a dense subspace of $\mathcal{H}$. If $\ell^p_m({\mathbb{Z}}^d)$ is not included in $\ell^2({\mathbb{Z}}^d)$ and $1 \leq p < \infty$ then we define $\mathcal{H}_m^p$ to be the completion of the subspace $\mathcal{H}_0$ of all finite linear combinations in $\cG$ with respect to the norm $\|f\|_{\mathcal{H}_m^p}= \| (\langle f,\tilde g_n \rangle )_{n \in {\mathbb{Z}}^d}\|_{\ell^p_m}$. If $p=\infty$ then we take the weak$^*$-completion of $\mathcal{H}_0$ to define $\mathcal{H}^\infty_m$.\\

\begin{rem} Under our assumptions one has $\mathcal{H}_m^p(\cG,\tilde{\cG}) = \mathcal{H}_m^p(\tilde \cG,\cG)$. Under the additional assumption that $\cA$ is a Banach $*$-algebra, the definition of $\mathcal{H}_m^p(\cG,\tilde{\cG})$ does even not depend on the particular $\cA$-self-localized dual chosen, and any other couple $(\cF,\tilde \cF)$ of $\cA$-self-localized dual frames which are localized to $\cG$ generates in fact the same spaces. See \cite{FoG,FR} for major details.
\end{rem}


Then, it is almost immediate to verify the following statement, see \cite{FoG}.

\begin{tm}
\label{banachfr}
Assume that $(\cG,\tilde \cG)$ is a pair of dual  $\cA$-self-localized frames for $\mathcal{H}$ with $\cG \sim_{\cA} \tilde \cG$. Then $\cG$ and its canonical dual frame $\tilde \cG$ are Banach frames and atomic decompositions for $\mathcal{H}_m^p(\cG,\tilde{\cG})$.
\end{tm}

\section{$\alpha$-modulation spaces}

\subsection{$\alpha$-modulation spaces as decomposition spaces}

In this section we want to recall the definition of $\alpha$-modulation spaces based on decomposition methods, without introducing them in full generality. For major details we refer to \cite{G,FG,F}.
In fact the spaces depend on a parameter $\alpha \in [0,1]$ which is a ``tuning tool'' to perform a suitable {\it segmentation} (decomposition) of the frequency domain as an {\it intermediate} geometry between those of modulation \cite{F4,Gr} and Besov \cite{FJ,Tr,Tr2} spaces.   

\begin{definition}
A countable set $\mathcal{I}$ of intervals $I \subset \mathbb{R}$ is called an \emph{admissible covering} of $\mathbb{R}$ if
\begin{itemize}
\item[(a)] $\mathbb{R} = \bigcup_{I \in \cI} I$, and
\item[(b)] $\#\{I \in \cI: x \in I\} \leq 2$ for all $x \in \mathbb{R}$.
\end{itemize}
Furthermore, if there exists a constant $0 \leq \alpha \leq 1$ such that $|I| \asymp (1+|\xi|)^\alpha$ for all $I \in \cI_\alpha$, and all $\xi \in I$, then $\cI_\alpha$ is called an $\alpha$-\emph{covering}.
\end{definition}
For an $\alpha$-covering $\mathcal{I}_\alpha$ one can identify the constituting intervals by means of two maps.
The \emph{position map} $p_\alpha$ from $\mathbb{Z}$ to $\mathbb{R}$, $p_\alpha: j \rightarrow p_\alpha(j)$, and the \emph{size map} $s_\alpha$ from $\mathbb{Z}$ to $\mathbb{R}_+$, $s_\alpha: j \rightarrow s_\alpha(j)$, so that the map from $\mathbb{Z}$ to $\cI_\alpha$, $j \rightarrow I_j$, $I_j=p_\alpha(j)+\text{sgn}(p_\alpha(j))[0,s_\alpha(j)]$ for $p_\alpha(j)\neq 0$, $I_j=[-s_\alpha(j),s_\alpha(j)]$ otherwise, is a bijection.
\begin{exmp}[Fornasier, Feichtinger \cite{FF}] 
For $b>0$ and $\alpha \in [0,1)$ an explicit example of $\alpha$-covering has been constructed in \cite{FF}, by choosing as position and size functions 
\begin{equation}
\label{pf}
p_\alpha(j) = \text{sgn}(j) \left( (1+(1-\alpha) \cdot b \cdot |j|  )^{\frac{1}{1-\alpha}} -1 \right )
\end{equation}
 and 
\begin{equation}
\label{sf}
s_\alpha(j) = b \cdot (1+(1-\alpha)\cdot b \cdot  (|j|+1))^{\frac{\alpha}{1-\alpha}},
\end{equation}
respectively. In particular, for $\alpha \rightarrow 1$ one has
$$
 \cI_1 = \{  \text{sgn}(j) \left ((e^{b |j|} -1)+[0,e^{b (|j|+1)}] \right )\}_{j \in \mathbb{Z}\backslash \{0\}} \cup \{[-e^b,e^b]\},
$$
is again an $\alpha$-covering, and for $b=\ln(2)$ is dyadic.
\end{exmp}

Without loss of generality we can assume that, associated to an admissible $\alpha$-covering $\cI_\alpha$, one can construct \cite[Theorem 4.2]{F} a corresponding \emph{bounded admissible partition of the unity} ({\tt BAPU}) $\Psi^\alpha=\{\psi_I^\alpha\}_{I \in \cI_\alpha}$ in $\mathcal{S}(\mathbb{R})$, i.e., 
\begin{itemize}
\item[(p1)] $\sup_{I \in \cI_\alpha} \|\psi_{I}^\alpha\|_{\mathcal{F} L^1} < \infty$,
\item[(p2)] $\text{supp}(\psi_{I}^\alpha) \subset I$ for all $I \in \cI_\alpha$, and
\item[(p3)] $\sum_{I \in \cI_\alpha} \psi_I^\alpha(\xi) =1 $ for all $\xi \in \mathbb{R}$.
\end{itemize}
\bigskip
Furthermore we define the {\it segmentation operator} $\cP_I^\alpha$ by
\begin{equation}
\label{seg}
\cP_I^\alpha (f) := \mathcal{F}^{-1} ( \psi_I^\alpha \mathcal{F} f), \quad I \in \cI_\alpha, \quad \text{for all }f \in \cS'(\mathbb{R}).
\end{equation}
In the following we will also write $\cP_j^\alpha:=\cP_{I_j}^\alpha$ and $\psi^\alpha_{j}:=\psi^\alpha_{I_j}$.
\begin{figure}[ht]
\hbox to \hsize {\hfill \epsfig{file=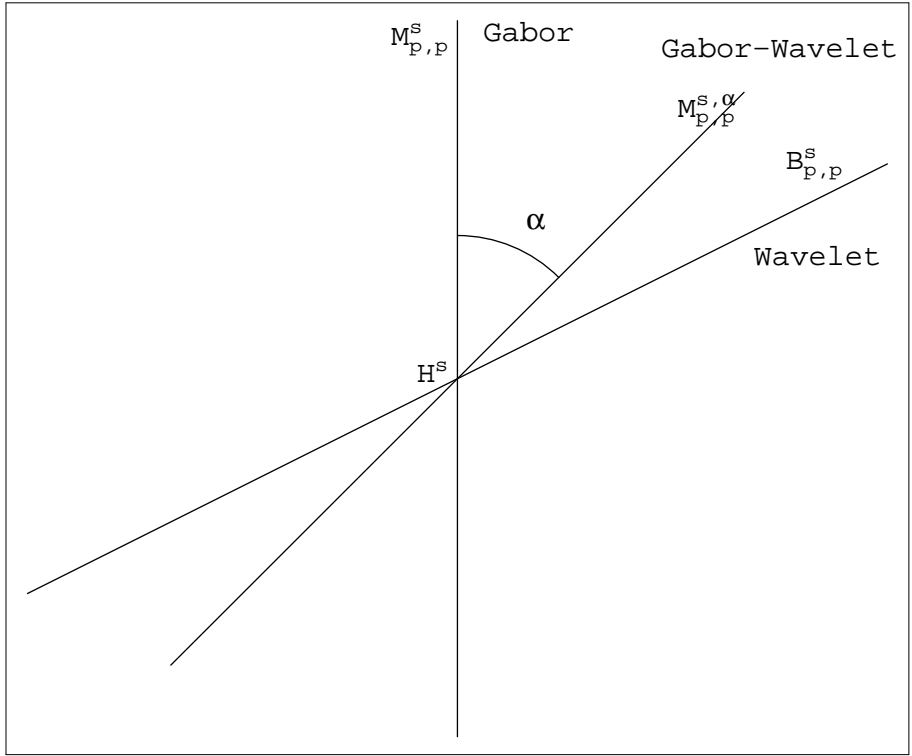, height=10cm} \hfill}
\caption{$\alpha$-modulation spaces}
\end{figure}

\begin{definition}[$\alpha$-modulation spaces, Gr\"obner \cite{G}]
\label{amod}
Given $1 \leq p,q \leq \infty$, $s \in \mathbb{R}$ and $0 \leq \alpha \leq 1$, let $\cI_\alpha$ be an $\alpha$-covering of $\mathbb{R}$ and let $\Psi^\alpha$ be a corresponding bounded admissible partition of the unity.
Then we define the $\alpha$-\emph{modulation space} $M_{p,q}^{s,\alpha}(\mathbb{R})$ for $q < \infty$ as the set of tempered distributions $f \in \cS'(\mathbb{R})$ satisfying 
\begin{equation}
\|f\|_{M_{p,q}^{s,\alpha}} := \L \sum_{I \in \cI_\alpha} \|\cP_I^\alpha (f)\|_p^q (1+|\omega_I|)^{s q} \R^{1/q} < \infty,
\end{equation}
with $\omega_I \in I$ for all $I \in \cI_\alpha$. For $q = \infty$ the definition is adapted substituting the $\ell^q$-norm with the $\sup$-norm over $I \in \cI_\alpha$. Let us denote $M_p^{s,\alpha}:= M_{p,p}^{s,\alpha}$.
\end{definition}

\begin{rem}
It is not difficult to check that the definition of $M_{p,q}^{s,\alpha}(\mathbb{R})$ does not depend on the particular choice of $\{\omega_I\}_{I \in \cI_\alpha}$. As a canonical choice we can assume $\omega_{I_j}=p_\alpha(j)$, for $I_j \in \cI_\alpha$. Moreover, two $\alpha$-coverings are equivalent in the sense of \cite[Definition 3.3]{FG}. A proof of such equivalence, even in higher dimension, can be found in \cite{G}.
As a consequence the definition of  $M_{p,q}^{s,\alpha}(\mathbb{R})$ does not depend on the particular choice of $\cI_\alpha$ \cite[Theorem 3.7]{FG} nor on $\{\cP_I^\alpha\}_{I \in \cI_\alpha}$ \cite[Theorem 2.3 (B)]{FG}.
In particular, from formula \eqref{pf}, we can assume without loss of generality that $p_\alpha(j) \asymp \text{sgn}(j) \left( (1+(1-\alpha) \cdot b \cdot |j|  )^{\frac{1}{1-\alpha}} -1 \right )$, $p_\alpha(0)=0$.
\end{rem}

\begin{exmps} 
\emph{Modulation spaces}. For $\alpha=0$ the space $M_{p,q}^{s,0}(\mathbb{R})$ coincides  with the {modulation space} $M_{p,q}^s(\mathbb{R})$. We refer to \cite{F4,Gr} for major details on such spaces.
They are naturally related to Gabor (time-frequency) frames, as we illustrate in the following.

 The combination of modulation and translation operators
\begin{equation}
\pi(\lambda) = M_{\omega} T_x \quad \mbox {for } \lambda=(x,\omega) \in \mathbb{R}^2
\end{equation}
is called a time-frequency shift. 
Let $\cX$ be a \emph{relatively separated} set in the time-frequency plane $\mathbb{R}^2$ and let $g \in L^2(\mathbb{R})$ be a fixed analyzing function. If the sequence $\mathcal{G}(g,\mathcal{X}) =\{\pi(\lambda) g\}_{\lambda \in \cX}$ is a frame for $L^2(\mathbb{R})$ then it is called {\it Gabor frame} if $\cX$ is a regular lattice, {\it non-uniform or irregular Gabor frame} otherwise.
If  $g \in  \cS(\mathbb{R})$ generates a (irregular) Gabor frame $\cG=\mathcal{G}(g,\mathcal{X})$ then for any $s>2$ the frame $\cG$ is intrinsically $s$-localized, i.e.,
$$
|\langle \pi(\lambda) g,  \pi(\mu) g \rangle| \lesssim (1+|\lambda-\mu|)^{-s}, \quad \lambda,\mu \in \cX,
$$
and, by \cite[Theorem 3.6, Corollary 3.7]{FG}, it has intrinsically $s$-self-localized canonical dual $\tilde{\cG}=\{\tilde e_\lambda\}_{\lambda \in \mathcal{X}}$. Moreover, it is shown in \cite{gro1,FoG} that $\cG$ and $\tilde{\cG}$ are Banach frames and atomic decompositions for suitable classes of modulation spaces. This means that
\begin{itemize} 
\item the frame expansions
\begin{equation}
f = \sum_{\lambda \in \mathcal{X}} \langle f, \tilde e_\lambda \rangle \pi(\lambda) g =  \sum_{\lambda \in \mathcal{X}} \langle f,  \pi(\lambda) g \rangle \tilde e_\lambda,
\end{equation}
converge unconditionally in $M_p^s(\mathbb{R})$;
\item the modulation space $M_p^s(\mathbb{R})$  can be characterized by the frame coefficients as follows:
\begin{equation}
\|f\|_{M_p^s} \asymp \|(\langle f, \tilde e_\lambda \rangle)_\lambda\|_{\ell^p_{m}(\mathcal{X})} \asymp \|(\langle f, \pi(\lambda) g)_\lambda\|_{\ell^p_{m}(\mathcal{X})}
\end{equation}
\end{itemize}
Therefore the spaces $\cH^p_{m}(\cG,\tilde \cG)$ and $M_p^s(\mathbb{R})$ coincide with equivalent norms, where here we have considered $m(\lambda)=m(x,\omega):=(1+|\omega|)^s$ as a polynomial weight depending only on the frequency variable.
\\

\noindent
\emph{Inhomogeneous Besov spaces}. For $\alpha \rightarrow 1$  the space $M_{p,q}^{s,1}(\mathbb{R})$ coincides with the {inhomogeneous Besov space} $B_{p,q}^{s}(\mathbb{R})$. Refer to \cite{FJ,Tr,Tr2} for major details on these classical spaces. It is well known \cite{M} that  inhomogeneous Besov spaces can be characterized by expansions of \emph{ wavelet frames} of the type $$\cG=\{T_k \varphi\}_{k \in \mathbb{Z}} \bigcup \{D_{2^{-j}} T_{k} \psi\}_{j \in \mathbb{N}, k \in \mathbb{Z}},$$
where $\varphi$ is a smooth refinable function and $\psi$ is a smooth wavelet function with enough vanishing moments. 

An application of the intrinsic localization of frame theory to characterize Besov space requires a different measure of localization. In particular, one should work with exponentially localized frames \cite{gro1,cg} as we will see also in the following. Therefore we postpone this limiting case to be discussed elsewhere. 
\end{exmps}

\subsection{Banach frames and atomic decompositions for $\alpha$-modulation spaces}
Assume $\alpha \in [0,1)$ and that $(p_\alpha,s_\alpha)$ is a pair of position and size functions.
Given the family
\begin{equation}
\cG:=\cG_\alpha(g,p_\alpha,s_\alpha,a)=\{M_{p_\alpha(j)} D_{s_\alpha^{-1}(j)} T_{a k} g\}_{j \in \mathbb{Z}, k \in \mathbb{Z}} \quad a>0,
\end{equation}
we want to illustrate under which (sufficient) conditions on the function $g$ one can ensure that $\cG$ is a frame for $L^2(\mathbb{R})$ and that $\cG$ extends also to a Banach frame and an atomic decomposition for a suitable family of Banach spaces. We want also to show that this class of Banach spaces is in fact constituted by $\alpha$-modulation spaces.
To this end, we discuss the properties of localization of $\cG$ and then we apply the principles illustrated in the previous section.
\\

\begin{rem}
For $\alpha=0$ the size function $s_0(j) \asymp (1+|p_0(j)|)^0=\text{const}$ and the position function $p_0$ describes a relatively separated set. Therefore, for $\alpha=0$ the frame $\cG$ is a Gabor frame. For $\alpha \rightarrow 1$, the dilation factor is controlled by $s_1(j)\asymp (1+|p_1(j)|)$. Therefore, since $\frac{p_1(j)}{s_1(j)} \asymp \text{const}$, the frame $\cG=\{e^{2 \pi i \frac{p_1(j)}{s_1(j)} a k} D_{s_1(j)^{-1}} T_{a k} (e^{2 \pi i \frac{p_1(j)}{s_1(j)} \cdot} g)\}_{j \in \mathbb{Z}, k \in \mathbb{Z}}$ is just a slight modification of a wavelet type frame.
\end{rem}
\\

Let us prove first some useful technical lemmas.

\begin{lemma} 
\label{lemma31}
Assume $s>1$.
\begin{itemize}
\item [(a)] For any  $0<\delta \leq 1$
\begin{equation}
\label{intdelta}
        \int_{\mathbb{R}} (1+|x-n|)^{-s} (\delta+|x-m|)^{-s} dx \lesssim \delta^{1-s}(\delta+|n-m|)^{-s}, \quad \text{for all } m,n \in \mathbb{R}.
\end{equation}
\item[(b)]  For $\rho\geq 1$ define $\Omega_\rho = \{x \in \mathbb{R}: |x| \geq  \rho\}$. For any $s'>\frac{1}{2}$ such that $s > s' + \frac{1}{2}$  and for any $b \geq 1$
$$
       \int_{\mathbb{R}}\left(  \chi_{\Omega_\rho} (b(x-n)) (1+|b(x-n)|)^{-s} (1+|x-m|)^{-s} \right) dx 
$$
\begin{equation}
\label{conv}
\leq C_\rho (1+|n-m|)^{-s'}, \quad \text{for all } m,n \in \mathbb{R},
\end{equation}
where $C_\rho \lesssim \left( \int_{\Omega_\rho} (1+|x|)^{-2(s-s')} dx \right )^{1/2}\rightarrow 0$, for $\rho \rightarrow 0$. In particular, $C_\rho \lesssim \rho^{1/2-(s-s')}$.
\end{itemize}
\end{lemma}
\begin{proof}
The statement (a) can be proved with similar arguments as \cite [Lemma 2.2]{gro1}:
Denote $A_1:=\{x \in \mathbb{R}:|n-x| \leq \frac{|n-m|}{2} \}$, and $A_2:=\mathbb{R} \backslash A_1$. If $x \in A_1$ then $|m-x| \geq \frac{|n-m|}{2}$ and
\begin{eqnarray*}
\int_{\mathbb{R}} (1+|x-n|)^{-s} (\delta+|x-m|)^{-s} dx &\leq& \left (\delta+\frac{|n-m|}{2} \right)^{-s} \int_{\mathbb{R}} (1+|x-n|)^{-s} dx\\
&\leq&2^{s} \left (\int_{\mathbb{R}} (1+|x|)^{-s} dx \right) (\delta+|n-m|)^{-s}.
\end{eqnarray*}
If $x \in A_2$ then $|n-x| > \frac{|n-m|}{2}$ and 
\begin{eqnarray*}
\int_{\mathbb{R}} (1+|x-n|)^{-s} (\delta+|x-m|)^{-s} dx &\leq& \left (1+\frac{|n-m|}{2} \right)^{-s} \int_{\mathbb{R}} (\delta+|x-m|)^{-s} dx\\
&\leq&2^{s} \delta^{-s} \left (\int_{\mathbb{R}} (1+|x/\delta|)^{-s} dx \right) (\delta+|n-m|)^{-s}\\
&=& 2^{s} \delta^{1-s} \left (\int_{\mathbb{R}} (1+|x|)^{-s} dx \right) (\delta+|n-m|)^{-s}.
\end{eqnarray*}
Therefore one has \eqref{intdelta}.

Let us prove (b). By assumption one has $\frac{w_{s'}}{w_s} \in L^2$, and this implies that $L^\infty_{w_s} \subset L^2_{w_{s'}}$, and  by Young inequality
\begin{equation}
\label{inc}
       L^2_{w_{s'}} \ast L^2_{w_{s'}} \subset L^\infty_{w_{s'}},
\end{equation}
where $\ast$ is the convolution operator.
The integral in \eqref{conv} can be interpreted as a convolution: Writing $w_{-s}^{\rho,b}(x) = \chi_{\Omega_\rho}(b x) w_{-s}(b x)$, one has
$$
 \left(w_{-s}^{\rho,b} \ast w_{-s}\right) (n-m) = \int_{\mathbb{R}}\left(  \chi_{\Omega_\rho} (b(x-n)) (1+|b(x-n)|)^{-s} (1+|x-m|)^{-s} \right) dx.
$$
By the continuous inclusion \eqref{inc}, a possible constant $C_\rho$ can be given by
$$
       C_\rho = \|w_{-s}^{\rho,b}\|_{L^2_{w_{s'}}} \|w_{-s}\|_{L^2_{w_{s'}}}. 
$$     
The norm 
\begin{eqnarray*}
\|w_{-s}^{\rho,b}\|_{L^2_{w_{s'}}} &=& \left (\int_{\mathbb{R}} \chi_{\Omega_\rho}(b x) (1+|bx|)^{-2s}(1+|x|)^{2 s'} dx \right)^{1/2}\\
&=& b^{-1/2} \left (\int_{\mathbb{R}} \chi_{\Omega_\rho}(x) (1+|x|)^{-2s}(1+|b^{-1}x|)^{2 s'} dx \right)^{1/2}\\
&\leq& \left (\int_{\Omega_\rho} (1+|x|)^{2(s'-s)} dx \right)^{1/2}
\end{eqnarray*}
Therefore $C_\rho \lesssim \rho^{1/2-(s-s')}$.
\end{proof}

\begin{rem}
Before proving the main technical lemma of this paper, it is useful to recall some properties of the fundamental operators of translation, modulation, and dilation, and of the pairs of position and size functions  $(p_\alpha,s_\alpha)$ we are going to consider.

1. With respect to the Fourier transform one has the following relations:
$$
\mathcal F M_\omega = T_{\omega} \mathcal F, \quad \mathcal F T_x = M_{-x} \mathcal F, \quad \mathcal F D_a = D_{a^{-1}} \mathcal F, \mbox{ for } x, \omega \in \mathbb{R}, a \in \mathbb{R}_+.
$$
One has also the following commutation relations:
$$
D_a T_x = T_{a x} D_a, \quad T_x M_\omega = e^{-2 \pi i x \omega} M_\omega T_x,  \mbox{ for } x,\omega \in \mathbb{R}, a \in \mathbb{R}_+.
$$

2. In the following we will assume that the pairs of position and size functions  $(p_\alpha,s_\alpha)$ we are going to consider satisfy the following properties for $|i| \leq |j|$, $i,j \in \mathbb{Z}$, 
\begin{itemize}
\item[(ps0)] $p_\alpha(j) j \geq 0$;
\item[(ps1)] $|p_\alpha(i)| \leq |p_\alpha(j)|$, $s_\alpha(i) \leq s_\alpha(j)$;
\item[(ps2)] $p_\alpha(i) / s_\alpha(i)= c(i)i$, for a suitable $(c(i) -(1-\alpha)) \asymp |i|^{-1}$ for $|i| \rightarrow \infty$;
\item[(ps3)] $\frac{|p_\alpha(j)| s_\alpha(j-\text{sgn}(j))}{|p_\alpha(j-\text{sgn}(j))| s_\alpha(j)} \geq 1$.
\end{itemize}
Of course, the position and size function in formulas \eqref{pf} and \eqref{sf} fulfill these requirements.
In particular for (ps3) it is sufficient to observe that for $x \in \mathbb R$ 
$$
        \lim_{x \rightarrow +\infty} \frac{p_\alpha(x) s_\alpha(x-1)}{p_\alpha(x-1) s_\alpha(x)} =1,
$$
$$
        \frac{p_\alpha(1) s_\alpha(0)}{p_\alpha(0) s_\alpha(1)} \geq 1,
$$
and that the derivative of $\frac{p_\alpha(x) s_\alpha(x-1)}{p_\alpha(x-1) s_\alpha(x)}$ with respect to $x$ is negative for $x \in [1, +\infty)$.
\end{rem}

\begin{lemma}
\label{lemma32}
Assume $0<a\leq 1$, $\gamma_f, \gamma_t > 1$, $\alpha \in [0,1)$, and let $(p_\alpha,s_\alpha)$ be a pair of position and size functions satisfying properties \emph{(ps0-3)}. \\

Let $\{g_\ell\}_{\ell \in \mathbb{Z}}, \{f_\ell\}_{\ell \in \mathbb{Z}} \subset L^1(\mathbb{R}) \cap C(\mathbb{R})$ such that
\begin{equation}
| g_\ell(x) \mathcal{F} g_\ell(\omega)|\lesssim (1+|x|)^{-\gamma_t}(1+|\omega|)^{-\gamma_f}, \quad x, \omega \in \mathbb{R},
\end{equation}
\begin{equation}
| f_\ell(x) \mathcal{F} f_\ell(\omega)|\lesssim (1+|x|)^{-\gamma_t}(1+|\omega|)^{-\gamma_f}, \quad x, \omega \in \mathbb{R},
\end{equation}
uniformly with respect to $\ell \in \mathbb{Z}$.
Then, 
\begin{itemize}  
\item[(a)] one has
$$
        |\langle  M_{ p_\alpha(j)}  D_{s_\alpha(j)^{-1}} T_{a k  } g_j,  M_{ p_\alpha(i)}  D_{s_\alpha(i)^{-1}} T_{a h  } f_i \rangle| 
$$
\begin{equation}
\label{loc1} 
\lesssim a^{-\frac{\gamma_t}{2}} \left (1 + |j-i| \right)^{\frac{1}{2}\left(\frac{\alpha}{(1-\alpha)}\gamma_t-\gamma_f\right)}   \left( 1+  \max\{s_\alpha(i),s_\alpha(j)\} |k s_\alpha(j)^{-1} - h s_\alpha(i)^{-1}| \right)^{-\frac{\gamma_t}{2}},
\end{equation}
for all  $i,j,h,k \in \mathbb{Z}$.
\item[(b)]  for a suitable system of segmentation operators $\{\mathcal{P}_j^\alpha\}_{j \in \mathbb{Z}}$ \eqref{seg} associated to a {\tt BAPU} $\Psi^\alpha=\{\psi_j^\alpha\}_{j \in \mathbb{Z}}$, one has
 $$
        |\langle  \mathcal{P}_j^\alpha M_{ p_\alpha(j)}  D_{s_\alpha(j)^{-1}} T_{a k  } g_j,  M_{ p_\alpha(i)}  D_{s_\alpha(i)^{-1}} T_{a h  } f_i \rangle|
$$
\begin{equation} 
\lesssim  a^{-\frac{\gamma_t}{2}} \left (1 + |j-i| \right)^{\frac{1}{2}\left(\frac{\alpha}{(1-\alpha)}\gamma_t-\gamma_f\right)}   \left( 1+  \max\{s_\alpha(i),s_\alpha(j)\} |k s_\alpha(j)^{-1} - h s_\alpha(i)^{-1}| \right)^{-\frac{\gamma_t}{2}},
\end{equation}
and
$$
        |\langle  \mathcal{P}_j^\alpha M_{ p_\alpha(j)}  D_{s_\alpha(j)^{-1}} T_{a k  } g_j,  \mathcal{P}_i^\alpha M_{ p_\alpha(i)}  D_{s_\alpha(i)^{-1}} T_{a h  } f_i \rangle| 
$$
\begin{equation} 
\lesssim  a^{-\frac{\gamma_t}{2}} \left (1 + |j-i| \right)^{\frac{1}{2}\left(\frac{\alpha}{(1-\alpha)}\gamma_t-\gamma_f\right)}   \left( 1+  \max\{s_\alpha(i),s_\alpha(j)\} |k s_\alpha(j)^{-1} - h s_\alpha(i)^{-1}| \right)^{-\frac{\gamma_t}{2}}.
\end{equation}
for all  $i,j,h,k \in \mathbb{Z}$.
\item[(c)] Let us consider $\rho \geq 1$ and  $\varphi \in C^\infty_c(\mathbb{R})$, $\text{supp}(\varphi)= [-(1+\varepsilon),1+\varepsilon]$,  with $\varphi \equiv 1$ on $[-1,1]$. Define $(g_\ell)_\rho:=\mathcal{F}^{-1}(\varphi(\frac{\cdot}{\rho}) \mathcal{F} g_\ell)$ a band-limited approximation of $g_\ell$ and $g_\ell^\rho := g_\ell - (g_\ell)_\rho$. For $\gamma_f' > 1$ and $\gamma_f > \gamma_f'+\gamma_t+3/2$, if $a=a(\rho) \asymp \rho^{-1}$, then
$$
       |\langle  M_{ p_\alpha(j)} T_{a \cdot s_\alpha(j)^{-1} \cdot k  } D_{s_\alpha(j)^{-1}} g_j^\rho,  M_{ p_\alpha(i)} T_{a \cdot s_\alpha(i)^{-1} \cdot h  } D_{s_\alpha(i)^{-1}} f_i \rangle| 
$$
\begin{equation} 
\lesssim D_\rho \left (1 + |j-i| \right)^{\frac{1}{2}\left(\frac{\alpha}{(1-\alpha)}\gamma_t-\gamma_f'\right)}   \left( 1+  \max\{s_\alpha(i),s_\alpha(j)\} |k s_\alpha(j)^{-1} - h s_\alpha(i)^{-1}| \right)^{-\frac{\gamma_t}{2}},
\end{equation}
for all  $i,j,h,k \in \mathbb{Z}$, where $D_\rho \rightarrow 0$ for $\rho \rightarrow \infty$, uniformly with respect to  $i,j,h,k \in \mathbb{Z}$.
\end{itemize}
\end{lemma}
\begin{proof}
Let us start showing (a), and, in particular, the case $j \geq i \geq 0$; the other cases can be shown with similar arguments. 
$$
|\langle  M_{ p_\alpha(j)} T_{a \cdot s_\alpha(j)^{-1} \cdot k  } D_{s_\alpha(j)^{-1}} g_j,  M_{ p_\alpha(i)} T_{a \cdot s_\alpha(i)^{-1} \cdot h  } D_{s_\alpha(i)^{-1}} f_i \rangle | 
$$
$$
= |\langle  M_{ -a k s_\alpha(j)^{-1}} T_{p_\alpha(j)} D_{s_\alpha(j)} \mathcal{F} g_j,   M_{ -a h s_\alpha(i)^{-1}} T_{p_\alpha(i)} D_{s_\alpha(i)} \mathcal{F} f_i \rangle| 
$$
\begin{equation}
\label{f17}             
        = \left |\int_{\mathbb{R}} \left( T_{p_\alpha(j)} D_{s_\alpha(j)} \mathcal{F} g_j(\omega) \right) \overline{ \left( T_{p_\alpha(i)} D_{s_\alpha(i)} \mathcal{F} f_i(\omega) \right)} e^{-2 \pi i \left(a (k s_\alpha(j)^{-1} - h s_\alpha(i)^{-1})\right) \omega} d\omega \right |.
\end{equation}
Step 1. (Frequency localization)\\

From \eqref{f17}  one has an estimation of \eqref{loc1} in the frequency domain:
$$
|\langle  M_{ p_\alpha(j)}  D_{s_\alpha(j)^{-1}} T_{a k  } g_j,  M_{ p_\alpha(i)}  D_{s_\alpha(i)^{-1}} T_{a h  } f_i \rangle|  \leq \int_{\mathbb{R}} \left| T_{p_\alpha(j)} D_{s_\alpha(j)} \mathcal{F} g_j(\omega) T_{p_\alpha(i)} D_{s_\alpha(i)} \mathcal{F} f_i(\omega) \right| d\omega 
$$
\begin{eqnarray*}
        &\lesssim&  \left(\frac{1}{s_\alpha(j) s_\alpha(i)}\right)^{1/2} \int_{\mathbb{R}} \L 1+\left |\frac{\omega - p_\alpha(j)}{s_\alpha(j)} \right| \R^{-\gamma_f}  \L1+\left |\frac{\omega - p_\alpha(i)}{s_\alpha(i)}\right|\R^{-\gamma_f} d\omega \\
&=&   \left(\frac{1}{s_\alpha(j) s_\alpha(i)}\right)^{1/2} \int_{\mathbb{R}} \left(1+\left |\frac{\omega}{s_\alpha(j)} - \frac{p_\alpha(j)}{s_\alpha(j)} \right| \right)^{-\gamma_f}  \\
&\times&\left (1+\left |\frac{\omega}{s_\alpha(j)}\frac{s_\alpha(j)}{s_\alpha(i)} - \frac{s_\alpha(j)}{s_\alpha(i)} \frac{s_\alpha(i)}{s_\alpha(j)} \frac{p_\alpha(i)}{s_\alpha(i)}\right|\right)^{-\gamma_f} d\omega\\
        &=&  \left ( \frac{s_\alpha(j)}{s_\alpha(i)}\right)^{1/2}\int_{\mathbb{R}} \left (1+\left |\omega- \frac{p_\alpha(j)}{s_\alpha(j)} \right| \right)^{-\gamma_f} \left (1+\frac{s_\alpha(j)}{s_\alpha(i)} \left |\omega - \frac{p_\alpha(i)}{s_\alpha(j)}\right| \right)^{-\gamma_f} d\omega 
\end{eqnarray*}
\begin{eqnarray}
\label{f18}
&\lesssim&   \left ( \frac{s_\alpha(j)}{s_\alpha(i)}\right)^{1/2}\int_{\mathbb{R}} \left (1+\left |\omega- \frac{p_\alpha(j)}{s_\alpha(j)} \right| \right)^{-\gamma_f} \left (1+\left |\omega - \frac{p_\alpha(i)}{s_\alpha(j)}\right| \right)^{-\gamma_f} d\omega.
\end{eqnarray}
By property (ps3) one has also that 
$$
\frac{p_\alpha(j)}{s_\alpha(j)} - \frac{p_\alpha(i)}{s_\alpha(i)} \geq 0.
$$
This implies, by property (ps1), the following inequality 
\begin{eqnarray*}
 \left |\frac{p_\alpha(j)}{s_\alpha(j)} - \frac{p_\alpha(i)}{s_\alpha(i)}\right| = \frac{p_\alpha(j)}{s_\alpha(j)} - \frac{p_\alpha(i)}{s_\alpha(i)}
\leq \frac{p_\alpha(j)}{s_\alpha(j)} - \frac{p_\alpha(i)}{s_\alpha(j)} 
= \left | \frac{p_\alpha(j)}{s_\alpha(j)} - \frac{p_\alpha(i)}{s_\alpha(j)} \right|.
\end{eqnarray*}
An application of Lemma \ref{lemma31} (a) and this last inequality  give
$$
\eqref{f18} \lesssim \left ( \frac{s_\alpha(j)}{s_\alpha(i)}\right)^{1/2} \left (1+\left |\frac{p_\alpha(j)}{s_\alpha(j)} - \frac{p_\alpha(i)}{s_\alpha(j)}\right| \right)^{-\gamma_f} \lesssim \left ( \frac{s_\alpha(j)}{s_\alpha(i)}\right)^{1/2} \left (1+\left |\frac{p_\alpha(j)}{s_\alpha(j)} - \frac{p_\alpha(i)}{s_\alpha(i)}\right| \right)^{-\gamma_f}. 
$$
Observing that $\frac{1+|y|}{1+|x|}\leq (1+|x-y|)$ for all $x,y \in \mathbb{R}$, one has by property (ps2)
\begin{eqnarray*}
        & &|\langle  M_{ p_\alpha(j)}  D_{s_\alpha(j)^{-1}} T_{a k  } g_j,  M_{ p_\alpha(i)}  D_{s_\alpha(i)^{-1}} T_{a h  } f_i \rangle|  \\
&\lesssim&  \left (1+|j-i| \right)^{\frac{\alpha}{2(1-\alpha)}} \left (1 + |j-i| \right)^{-\gamma_f} 
\end{eqnarray*}
\begin{eqnarray}
\label{f21}
&=& \left (1 + |j-i| \right)^{\frac{\alpha}{2(1-\alpha)}-\gamma_f}.
\end{eqnarray}
Step 2. (Time localization)\\

From \eqref{f17} one has an estimation of \eqref{loc1} also in the time domain: 
\begin{eqnarray*}
&&|\langle  M_{ p_\alpha(j)}  D_{s_\alpha(j)^{-1}} T_{a k  } g_j,  M_{ p_\alpha(i)}  D_{s_\alpha(i)^{-1}} T_{a h  } f_i \rangle| \\
&\leq & (D_{s_\alpha(i)^{-1}}|f_i|) \ast(D_{s_\alpha(j)^{-1}} |g_j|)^\nabla(a (k s_\alpha(j)^{-1} - h s_\alpha(i)^{-1})) 
\end{eqnarray*}
\begin{eqnarray}
\label{f22}
&\lesssim&  \left (s_\alpha(j) s_\alpha(i)\right)^{1/2} \int_\mathbb{R} \left( 1+ |s_\alpha(j)( y -x)|\right)^{-\gamma_t} \left( 1 + |s_\alpha(i)x|\right)^{-\gamma_t} dx,
\end{eqnarray}
where $y = a (k s_\alpha(j)^{-1} - h s_\alpha(i)^{-1})$.
By a change of variable and Lemma \ref{lemma31} (a), formula \eqref{f22} can be expressed and then estimated by
\begin{eqnarray*}
&& \left ( \frac{s_\alpha(j)}{s_\alpha(i)}\right)^{1/2} \int_\mathbb{R} \left( 1+ |s_\alpha(j)(y -\frac{x}{s_\alpha(i)}|\right)^{-\gamma_t} \left( 1 + |x|\right)^{-\gamma_t} dx\\
&=&  \left ( \frac{s_\alpha(j)}{s_\alpha(i)}\right)^{1/2-\gamma_t} \int_\mathbb{R} \left( \frac{s_\alpha(i)}{s_\alpha(j)} + |s_\alpha(i) y -x|\right)^{-\gamma_t} \left( 1 + |x|\right)^{-\gamma_t} dx \\
       &\lesssim& \left (1+|j-i| \right)^{\frac{\alpha}{2(1-\alpha)}} \left ( \frac{s_\alpha(j)}{s_\alpha(i)}\right)^{\gamma_t-1} \left ( \frac{s_\alpha(j)}{s_\alpha(i)}\right)^{-\gamma_t} \left( \frac{s_\alpha(i)}{s_\alpha(j)} + |s_\alpha(i) y| \right)^{-\gamma_t}
\end{eqnarray*}
\begin{eqnarray}
\label{f23}
 &\lesssim& \left (1+|j-i| \right)^{\frac{\alpha}{(1-\alpha)}(\gamma_t-1/2)}  \left(1 + |s_\alpha(j) y| \right)^{-\gamma_t}.
\end{eqnarray}
Step 3. (Time-frequency localization)\\
By combining formulae \eqref{f21} and \eqref{f23}, and assuming $a\leq 1$, one has
\begin{eqnarray*}
 &&|\langle  M_{ p_\alpha(j)} T_{a \cdot s_\alpha(j)^{-1} \cdot k  } D_{s_\alpha(j)^{-1}} g_j,  M_{ p_\alpha(i)} T_{a \cdot s_\alpha(i)^{-1} \cdot h  } D_{s_\alpha(i)^{-1}} f_i \rangle |^2  \\
&\lesssim&  \left (1 + |j-i| \right)^{\frac{\alpha}{(1-\alpha)}\gamma_t-\gamma_f}   \left( 1+  a \max\{s_\alpha(i),s_\alpha(j)\} |k s_\alpha(j)^{-1} - h s_\alpha(i)^{-1}| \right)^{-\gamma_t} \\
&\lesssim&  a^{-\gamma_t} \left (1 + |j-i| \right)^{\frac{\alpha}{(1-\alpha)}\gamma_t-\gamma_f}   \left( 1+  \max\{s_\alpha(i),s_\alpha(j)\} |k s_\alpha(j)^{-1} - h s_\alpha(i)^{-1}| \right)^{-\gamma_t}
\end{eqnarray*}
We want to show now (b).\\

Observe that
$$
         \mathcal{F}(\mathcal{P}_j^\alpha M_{ p_\alpha(j)} T_{a \cdot s_\alpha(j)^{-1} \cdot k  } D_{s_\alpha(j)^{-1}} g_j) =  \psi_j^\alpha T_{p_\alpha(j)}  M_{-a k s_\alpha(j)^{-1}} D_{s_\alpha(j)} \mathcal{F} g_j,
$$
Without loss of generality, by similar arguments as in \cite[Theorem 4.2]{F} we can assume $\psi_j^\alpha = s_\alpha(j)^{1/2} T_{p_\alpha(j)} D_{s_\alpha(j)} \varphi_j^\alpha$, with
$$
\varphi_j^\alpha(x) \mathcal{F} \varphi_j^\alpha(\omega) \lesssim(1+|x|)^{-\gamma_t}(1+|\omega|)^{-\gamma_f}
$$
for all $j \in \mathbb{Z}$ and $x, \omega \in \mathbb{R}$.
Therefore 
$$
        \mathcal{F}(\mathcal{P}_j^\alpha M_{ p_\alpha(j)} T_{a \cdot s_\alpha(j)^{-1} \cdot k  } D_{s_\alpha(j)^{-1}} g_j) =   T_{p_\alpha(j)} M_{- a k s_\alpha(j)^{-1}} D_{s_\alpha(j)} (\varphi_j^\alpha \mathcal{F} g_j).
$$
If $|\mathcal{F} g_j(\omega)| \lesssim (1+|\omega|)^{-\gamma_f}$, then $|\varphi_j^\alpha \mathcal{F} g_j(\omega)| \lesssim (1+|\omega|)^{-\gamma_f}$, uniformly with respect to $j \in \mathbb{Z}$.
Moreover, by Lemma \ref{lemma31} (a), one has
$$
        |\mathcal{F}^{-1} \L \varphi_j^\alpha \mathcal{F} g_j \R (x)| \lesssim (1+|x|)^{-\gamma_t},
$$
uniformly with respect to $j \in \mathbb{Z}$.
At this point one can conclude the proof of (b) by an application of (a).
\\

Let us now show the last statement (c). First of all observe that
\begin{eqnarray*}
| g_\ell^\rho(x)| &\leq& |g_\ell(x)| + | (g_\ell)_\rho(x)| \\
&\leq& |g_\ell(x)| + |g_\ell\ast (\mathcal{F}^{-1} \varphi(\frac{\cdot}{\rho}))(x)|,
\end{eqnarray*}
and
\begin{eqnarray*}
\mathcal{F}^{-1} \varphi(\frac{\cdot}{\rho})(x)&=&\left |\int_{-\rho(1+\varepsilon)}^{\rho(1+\varepsilon)}  \varphi(\frac{\omega}{\rho}) e^{2 \pi i \omega x} d\omega \right|\\
&=& \rho \left |\int_{-(1+\varepsilon)}^{(1+\varepsilon)}  \varphi(\omega) e^{2 \pi i \omega x \rho} d\omega \right| \\
&\lesssim& \rho (1+|\rho x|)^{-\gamma_t} \leq  \rho (1+|x|)^{-\gamma_t} 
\end{eqnarray*}
By combining these two estimations and applying Lemma \ref{lemma31} (a) one has
\begin{equation}
\label{f25}
| g_\ell^\rho(x)| \lesssim \rho  (1+|x|)^{-\gamma_t}.
\end{equation}
Moreover, one has also the following estimation in the frequency
\begin{equation}
\label{f26}
|\mathcal{F} g_\ell^\rho(\omega)| = |\mathcal{F} g_\ell (\omega) (1-\varphi(\frac{\omega}{\rho}))| \lesssim \chi_{\Omega_\rho}(\omega)(1+|\omega|)^{-\gamma_f}.
\end{equation}
Estimations \eqref{f25} and \eqref{f26} yield
\begin{equation}
\label{f27}
| g_\ell^\rho(x) \mathcal{F} g_\ell^\rho(\omega)|\lesssim \rho (1+|x|)^{-\gamma_t}\chi_{\Omega_\rho}(\omega)(1+|\omega|)^{-\gamma_f}, \quad x, \omega \in \mathbb{R}.
\end{equation}
Following the computations done for the statement (a) in Step 1, one obtains one of the following two expressions depending, respectively, on the assumption, e.g., that $0\leq i\leq j$ or $0\leq j\leq i$.

\noindent Or
$$
|\langle  M_{ p_\alpha(j)}  D_{s_\alpha(j)^{-1}} T_{a k  } g_j^\rho,  M_{ p_\alpha(i)}  D_{s_\alpha(i)^{-1}} T_{a h  } f_i \rangle|  
$$
$$
\leq \left ( \frac{s_\alpha(j)}{s_\alpha(i)}\right)^{1/2}\int_{\mathbb{R}} \chi_{\Omega_\rho} \left(\omega- \frac{p_\alpha(j)}{s_\alpha(j)}\right ) \left(1+\left |\omega- \frac{p_\alpha(j)}{s_\alpha(j)} \right| \right)^{-\gamma_f} \left (1+\frac{s_\alpha(j)}{s_\alpha(i)} \left |\omega - \frac{p_\alpha(i)}{s_\alpha(j)}\right| \right)^{-\gamma_f} d\omega,
$$
Or
$$
|\langle  M_{ p_\alpha(j)}  D_{s_\alpha(j)^{-1}} T_{a k  } g_j^\rho,  M_{ p_\alpha(i)}  D_{s_\alpha(i)^{-1}} T_{a h  } f_i \rangle|  
$$
$$
\leq \left ( \frac{s_\alpha(i)}{s_\alpha(j)}\right)^{1/2}\int_{\mathbb{R}}   \chi_{\Omega_\rho} \left(\frac{s_\alpha(i)}{s_\alpha(j)} \left (\omega - \frac{p_\alpha(j)}{s_\alpha(i)}\right) \right )\left (1+\frac{s_\alpha(i)}{s_\alpha(j)} \left |\omega - \frac{p_\alpha(j)}{s_\alpha(i)}\right| \right)^{-\gamma_f}  \left(1+\left |\omega- \frac{p_\alpha(i)}{s_\alpha(i)} \right| \right)^{-\gamma_f}d\omega.
$$
In both the cases one can apply Lemma \ref{lemma31} (b) and conclude, as in Step 1, that
$$
|\langle  M_{ p_\alpha(j)}  D_{s_\alpha(j)^{-1}} T_{a k  } g_j^\rho,  M_{ p_\alpha(i)}  D_{s_\alpha(i)^{-1}} T_{a h  } f_i \rangle|  \lesssim C_\rho \left (1 + |j-i| \right)^{\frac{\alpha}{2(1-\alpha)}-\gamma_f'},
$$
where $C_\rho \lesssim \rho^{1/2-(\gamma_f-\gamma_f')}$.
Moreover, proceeding as in Step 2, and using the estimation \eqref{f25}, one obtains
$$
|\langle  M_{ p_\alpha(j)}  D_{s_\alpha(j)^{-1}} T_{a k  } g_j^\rho,  M_{ p_\alpha(i)}  D_{s_\alpha(i)^{-1}} T_{a h  } f_i \rangle| \lesssim \rho \left (1+|j-i| \right)^{\frac{\alpha}{(1-\alpha)}(\gamma_t-1/2)}  \left(1 + |\max\{s_\alpha(i),s_\alpha(j)\} y| \right)^{-\gamma_t}.
$$
Again, combining the last expressions one has
$$
|\langle  M_{ p_\alpha(j)}  D_{s_\alpha(j)^{-1}} T_{a k  } g_j^\rho,  M_{ p_\alpha(i)}  D_{s_\alpha(i)^{-1}} T_{a h  } f_i \rangle|^2
$$
$$
 \lesssim \rho C_\rho a^{-\gamma_t} \left (1 + |j-i| \right)^{\frac{\alpha}{(1-\alpha)}\gamma_t-{\gamma_f'}}   \left( 1+  \max\{s_\alpha(i),s_\alpha(j)\} |k s_\alpha(j)^{-1} - h s_\alpha(i)^{-1}| \right)^{-\gamma_t}
$$
Since we assume $a=a(\rho) \asymp \rho^{-1}$, one finally has
$$
|\langle  M_{ p_\alpha(j)}  D_{s_\alpha(j)^{-1}} T_{a k  } g_j^\rho,  M_{ p_\alpha(i)}  D_{s_\alpha(i)^{-1}} T_{a h  } f_i \rangle|^2
$$
$$
 \lesssim \rho^{3/2+\gamma_t-(\gamma_f-\gamma_f')}\left (1 + |j-i| \right)^{\frac{\alpha}{(1-\alpha)}\gamma_t-{\gamma_f'}}   \left( 1+  \max\{s_\alpha(i),s_\alpha(j)\} |k s_\alpha(j)^{-1} - h s_\alpha(i)^{-1}| \right)^{-\gamma_t},
$$
and $D_\rho^2:=\rho^{3/2+\gamma_t-(\gamma_f-\gamma_f')} \rightarrow 0$ for $\rho \rightarrow +\infty$.
\end{proof}


Inspired by the results of the previous technical lemma we state the following definition.

\begin{definition}
\label{locclass}
For $\alpha\in [0,1)$, $\gamma,\eta>1$ we define the class of the $(\alpha,\gamma,\eta)$-off-diagonal-decaying matrices $\cA_{\alpha,\gamma,\eta}$ on $\mathbb{Z}^2\times \mathbb{Z}^2$ as follows.
A matrix $A=(a_{jk,ih})_{i,j,h,k \in \mathbb{Z}} \in \cA_{\alpha,\gamma,\eta}$ if and only if
$$
| a_{jk,ih}| \leq K \left (1 + (1-\alpha)|j-i| \right)^{-\frac{\gamma}{1-\alpha}}   \left( 1+  \max\{s_\alpha(i),s_\alpha(j)\} |k s_\alpha(j)^{-1} - h s_\alpha(i)^{-1}| \right)^{-\eta},
$$
for a suitable $K>0$ constant independent on $i,j,h,k \in \mathbb{Z}$.
\end{definition}

\begin{rem}
Observe that this class of matrices is for $\alpha=0$ a Banach $*$-algebra, see \cite{j,gro1,GL}, typically arising in the localization theory of Gabor frames, see, e.g., Examples 1 and \cite{FoG}. It is also known that, for the case $\alpha \rightarrow 1$, i.e., the matrices localized as follows
$$
| a_{jk,ih}| \leq K e^{-\gamma |j-i|}   \left( 1+  \max\{e^i,e^j\} |k e^{-j} - h e^{-i}| \right)^{-\eta},
$$
cannot form an algebra, see for example \cite{Le,DFR}. This class of matrices typically arises in the localization theory of wavelet frames. In general, for $\alpha \in (0,1)$ it is not yet known whether $\cA_{\alpha,\gamma,\eta}$ can be an algebra. Interesting related results can be found in \cite{HN}.
\end{rem}

Of course, in order to use the localization concept for the characterization of Banach spaces, we should show that the matrices belonging to the class $\cA_{\alpha,\gamma,\eta}$ can be bounded on suitable weighted $\ell^p(\mathbb{Z}^2)$ spaces.
\begin{prop}
\label{boundedness}
Let $\alpha \in [0,1), \gamma,\eta>1$ be fixed. Then any  matrix $A \in \cA_{\alpha,(1-\alpha)\gamma,\eta}$ extends to a bounded operator from $\ell^p_m(\mathbb{Z}^2)$ to $\ell^p_m(\mathbb{Z}^2)$ for all $p \in [1,\infty]$ and for any $s$-moderate weight $m(j,k):=m(j)$, depending only on the first index, $0 \leq s<\gamma-1$. Moreover, one can estimate the operator norm  by $\|A\|_{\ell^p_m \rightarrow \ell^p_m} \lesssim K$, where $K$ is the constant appearing in the Definition \ref{locclass}.
\end{prop}
\begin{proof}
We first show that $A$ is bounded on  $\ell^1_m(\mathbb{Z}^2)$ and on  $\ell^\infty_m(\mathbb{Z}^2)$, and then we conclude by interpolation the boundedness on $\ell^p_m(\mathbb{Z}^2)$.
Consider $c \in \ell^1_m(\mathbb{Z}^2)$.
\begin{eqnarray*}
&& \|A c \|_{\ell^1_m(\mathbb{Z}^2)}\\
&\lesssim& K \sum_{j,k \in \mathbb{Z}} \left(\sum_{i,h \in \mathbb{Z}}  \left (1 + |j-i| \right)^{-\gamma}   \left( 1+  \max\{s_\alpha(i),s_\alpha(j)\} |k s_\alpha(j)^{-1} - h s_\alpha(i)^{-1}| \right)^{-\eta} |c_{i,h}| \right )m(j)\\
&=& K \sum_j \sum_i \left (1 + |j-i| \right)^{-\gamma} m(j) \left(\sum_h \left (\sum_k \left( 1+  \max\{s_\alpha(i),s_\alpha(j)\} |k s_\alpha(j)^{-1} - h s_\alpha(i)^{-1}| \right)^{-\eta} \right ) |c_{i,h}|\right)\\
&\lesssim& K \sum_j \sum_i \left (1 + |j-i| \right)^{-\gamma} m(j) \left(\sum_h \left (\int_{\mathbb{R}} \left( 1+  \max\{s_\alpha(i),s_\alpha(j)\} |x s_\alpha(j)^{-1} - h s_\alpha(i)^{-1}| \right)^{-\eta} dx \right ) |c_{i,h}|\right)\\
&\lesssim& K \sum_j \sum_i \left (1 + |j-i| \right)^{-\gamma} m(j) 
\left(\sum_h \left (\int_{\mathbb{R}} \left( 1+  |x| \right)^{-\eta} dx \right ) |c_{i,h}|\right)\\
&\lesssim& K \sum_j \sum_i \left (1 + |j-i| \right)^{-\gamma}\left(\sum_h  |c_{i,h}|\right) m(j).
\end{eqnarray*}
Let us denote $d_{i}:=\left(\sum_h  |c_{i,h}|\right)$. Of course $d=(d_i)_{i \in \mathbb{Z}} \in \ell^1_m(\mathbb{Z})$, and by \cite[Lemma 2.3]{gro1}
\begin{eqnarray*}
 \|A c \|_{\ell^1_m}&\lesssim& K \sum_j \left(\sum_i \left (1 + |j-i| \right)^{-\gamma} d_i \right) m(j)\\
&\lesssim& K \sum_j d_j m(j) = K \|c\|_{\ell^1_m(\mathbb{Z}^2)}.
\end{eqnarray*}
Similarly one can show the boundedness on $\ell^\infty_m(\mathbb{Z}^2)$. Consider $c \in \ell^\infty_m(\mathbb{Z}^2)$.
\begin{eqnarray*}
&& \|A c \|_{\ell^\infty_m(\mathbb{Z}^2)}\\
&\lesssim& K \sup_{j,k \in \mathbb{Z}} \left(\sum_{i,h \in \mathbb{Z}}  \left (1 + |j-i| \right)^{-\gamma}   \left( 1+  \max\{s_\alpha(i),s_\alpha(j)\} |k s_\alpha(j)^{-1} - h s_\alpha(i)^{-1}| \right)^{-\eta} |c_{i,h}| \right )m(j)\\
&\leq& K \sup_j \sum_i \left (1 + |j-i| \right)^{-\gamma} m(j) \left(\sup_k \sum_h \left( 1+  \max\{s_\alpha(i),s_\alpha(j)\} |k s_\alpha(j)^{-1} - h s_\alpha(i)^{-1}| \right)^{-\eta} |c_{i,h}|\right).
\end{eqnarray*}
Since we have already shown that
$$
 \sup_k \sum_h \left( 1+  \max\{s_\alpha(i),s_\alpha(j)\} |k s_\alpha(j)^{-1} - h s_\alpha(i)^{-1}| \right)^{-\eta} \lesssim 1,
$$
then 
$$
\|A c \|_{\ell^\infty_m(\mathbb{Z}^2)} \lesssim K \sup_j \sum_i (1 + |j-i|)^{-\gamma} m(j) \left( \sup_h |c_{i,h}|\right).
$$
Again, let us denote $d_{i}:=\left(\sup_h  |c_{i,h}|\right)$. Of course $d=(d_i)_{i \in \mathbb{Z}} \in \ell^\infty_m(\mathbb{Z})$, and  by \cite[Lemma 2.3]{gro1}
\begin{eqnarray*}
 \|A c \|_{\ell^\infty_m}&\lesssim& K \sup_j \left(\sum_i \left (1 + |j-i| \right)^{-\gamma} d_i \right) m(j)\\
&\lesssim& K \sup_j d_j m(j) = K \|c\|_{\ell^\infty_m(\mathbb{Z}^2)}.
\end{eqnarray*}
One concludes the proof by interpolation of $\ell^p_m(\mathbb{Z}^2)$ spaces \cite{BL}.
\end{proof}

Finally, we have developed all the technical tools in order to show the main result of Banach frame and atomic decomposition for $\alpha$-modulation spaces as follows.
\\

Assume $s>0$ and $\alpha \in [0,1)$. We say that $g \in L^1(\mathbb{R}) \cap C(\mathbb{R})$ is {\it $(s;\alpha)$-localized}, if, for some $\gamma_f' > 2\left (1 + \frac{s}{1-\alpha} \right) + \frac{\alpha}{1-\alpha}\gamma_t$, $\gamma_t>2$, and $\gamma_f > \gamma_f'+\gamma_t+3/2$,
\begin{equation}
| g(x) \mathcal{F} g(\omega)|\lesssim (1+|x|)^{-\gamma_t}(1+|\omega|)^{-\gamma_f}, \quad x, \omega \in \mathbb{R}.
\end{equation}
Of course, Schwartz functions are $(s;\alpha)$-localized for all $s \geq 0$ and all $\alpha \in [0,1)$.
\begin{tm}
\label{chamod}
Let $\alpha \in [0,1)$, $s \in \mathbb{R}$. Assume that $g \in L^1(\mathbb{R}) \cap C(\mathbb{R})$ is $(|s|;\alpha)$-localized, $\mathcal{F} g(\omega) \neq 0$ for $\omega \in \Omega_0=[-1,1]$ and that $(p_\alpha,s_\alpha)$ is a pair of position and size functions satisfying conditions (ps0-3). Then, there exists $0<a_0\leq 1$ small enough such that for all $0<a\leq a_0$ the family
\begin{equation}
\cG:=\cG_\alpha(g,p_\alpha,s_\alpha,a)=\{M_{p_\alpha(j)} D_{s_\alpha^{-1}(j)} T_{a k} g\}_{j \in \mathbb{Z}, k \in \mathbb{Z}},
\end{equation}
\begin{itemize}
\item[(a)] is  a $\cA_{\alpha,\frac{1-\alpha}{2} (\gamma_f-\frac{\alpha}{1-\alpha}\gamma_t),\frac{\gamma_t}{2}}$-self-localized frame for $L^2(\mathbb{R})$;
\item[(b)] is an atomic decomposition for the $\alpha$-modulation space $M^{s+\alpha(1/p-1/2)}_p$ for all $p \in [1,\infty]$;
\item[(c)] is a Banach frame for the $\alpha$-modulation space $M^{s+\alpha(1/p-1/2)}_p$ for all $p \in [1,\infty]$.
\end{itemize}
\end{tm}

\begin{proof}
The statement (a) is a direct consequence of an application of  \cite[Theorem 1]{FF} and Lemma \ref{lemma32} (a). Let us show (b,c).
\\

The proof is develop as follows: First we show that for a band-limited approximation $g_\rho$ of $g$ the system $\cG_\alpha(g_\rho,p_\alpha,s_\alpha,a)$ forms a Banach frame and an atomic decomposition for $M^{s+\alpha(1/p-1/2)}_p$, and then we extend the result to $\cG_\alpha(g,p_\alpha,s_\alpha,a)$ by the application of the perturbation results  \cite[Theorem 2.2, Theorem 2.3]{CH}.
\\

Denote $\cA :=\cA_{\alpha,\frac{1-\alpha}{2} (\gamma_f-\frac{\alpha}{1-\alpha}\gamma_t),\frac{\gamma_t}{2}}$.  
Let us consider $\rho \geq 1$ and  $\varphi \in C^\infty_c(\mathbb{R})$, $\text{supp}(\varphi)= [-(1+\varepsilon),1+\varepsilon]$,  with $\varphi \equiv 1$ on $[-1,1]$. Define $g_\rho:=\mathcal{F}^{-1}(\varphi(\frac{\cdot}{\rho}) \mathcal{F} g)$ a band-limited approximation of $g$ and $g^\rho := g - g_\rho$.
If $f \in M^{s+\alpha(1/p-1/2),\alpha}_{p}(\mathbb{R})$ then, for $j \in \mathbb{Z}$, $\mathcal{P}_j^\alpha(f)$ is an $L^p(\mathbb{R})$ band-limited function and, by classical theorems on series expansions of band-limited functions (see also \cite{FG5},\cite[Example 5]{Fo1})), there exists $a=a(\rho)\asymp \rho^{-1}$ such that
\begin{equation}
        \mathcal{P}_j^\alpha(f) = \sum_{k \in \mathbb{Z}} \langle  \mathcal{P}_j^\alpha f, M_{p_\alpha(j)} D_{s_\alpha^{-1}(j)} T_{a k} \tilde g_\rho \rangle M_{p_\alpha(j)} D_{s_\alpha^{-1}(j)} T_{a k} g_\rho,
\end{equation}
where $ \tilde g_\rho = a \tilde g$ is a well-decaying  band-limited dual function with $\mathcal{F} \tilde g \mathcal{F} g \equiv 1$ on $\Omega_0$, and
\begin{equation}
\label{equivbandnorm}
        s_\alpha(j)^{\frac{2-p}{2}} \cdot \|\mathcal{P}_j^\alpha(f)\|_p^p \asymp \sum_{k \in \mathbb{Z}} |\langle f,  \mathcal{P}_j^\alpha M_{p_\alpha(j)} D_{s_\alpha^{-1}(j)} T_{a k} \tilde g_\rho\rangle|^p,
\end{equation}
for $p<\infty$ and similarly one has the equivalence for $p=\infty$.
In particular, since $\mathcal{P}_j^\alpha f$ is band-limited and recalling that $a=a(\rho) \asymp \rho^{-1}$, one has
\begin{eqnarray*}
\sum_{k \in \mathbb{Z}} |\langle f,  \mathcal{P}_j^\alpha M_{p_\alpha(j)} D_{s_\alpha^{-1}(j)} T_{a k} \tilde g_\rho\rangle|^p &=& \sum_{k \in \mathbb{Z}} |\langle    \left(D_{s_\alpha(j)} M_{-p_\alpha(j)}\mathcal{P}_j^\alpha\right) f,   T_{a k} \tilde g_\rho\rangle|^p\\
&=&  a^p \sum_{k \in \mathbb{Z}} |  \left(D_{s_\alpha(j)} M_{-p_\alpha(j)}\mathcal{P}_j^\alpha\right) f \ast \tilde g^\nabla (a k)|^p\\
&\lesssim& a^{p-1} \|\left(D_{s_\alpha(j)} M_{-p_\alpha(j)}\mathcal{P}_j^\alpha\right) f \ast \tilde g^\nabla\|_p^p\\
&\lesssim&    s_\alpha(j)^{\frac{2-p}{2}} \cdot \|\mathcal{P}_j^\alpha(f)\|_p^p,
\end{eqnarray*}
uniformly with respect to $\rho \geq 1$ (see also \cite{FG5},\cite[Theorem 4, Remark 2]{FF},\cite[Example 5]{Fo1}).
Here we have used the fact that for an $L^p$-band-limited function $h$, $\|(h(a k))_{k\in \mathbb{Z}}\|_{\ell^p} \leq C a^{-1/p} \|h\|_p$.
The usual modifications apply for the case $p=\infty$.
By an application of \cite[Theorem 1]{FF} or \cite[Theorem 14 and Corollary 17]{Fo1}, the systems  
\begin{equation}
\cG_\rho:=\{M_{p_\alpha(j)} D_{s_\alpha^{-1}(j)} T_{a k} g_\rho\}_{j \in \mathbb{Z}, k \in \mathbb{Z}} \text{ and }\tilde \cG_\rho:=\{       \mathcal{P}_j^\alpha M_{p_\alpha(j)} D_{s_\alpha^{-1}(j)} T_{a k} \tilde g_\rho\}_{j \in \mathbb{Z}, k \in \mathbb{Z}}
\end{equation} 
constitute a dual pair $(\cG_\rho,\tilde \cG_\rho)$ of frames for $L^2(\mathbb{R})$. By Lemma \ref{lemma32} (a),(b)  $\cG_\rho$ and $\tilde \cG_\rho$ are $\cA$-self-localized and $\cG_\rho \sim_{\cA} \tilde \cG_\rho$. Therefore, by Proposition \ref{boundedness}, it makes sense to define the abstract Banach space $\cH^p_{m_{s,\alpha}}(\cG_\rho,\tilde \cG_\rho)$, where $m_{s,\alpha}(j,k):=m_{s,\alpha}(j)=(1+(1-\alpha)|j|)^\frac{s}{1-\alpha}$. By Definition \ref{amod} and \eqref{equivbandnorm}, one has the equivalence of norms:
\begin{equation}
\label{f33}
        \|f\|_{M^{ s + \alpha(\frac{1}{p}-\frac{1}{2}),\alpha}_{p}} \asymp \|\L\langle f, \mathcal{P}_j^\alpha M_{p_\alpha(j)} D_{s_\alpha^{-1}(j)} T_{a k} \tilde g_\rho\rangle \R_{j,k \in \mathbb{Z}} \|_{\ell^p_{m_{s,\alpha}}(\mathbb{Z}^2)} = \| f \|_{\cH^p_{m_{s,\alpha}}(\cG_\rho,\tilde \cG_\rho)}.
\end{equation}
It is not difficult to see that the space of linear combinations of elements of $\cG_\rho$ is in fact dense in $M^{ s + \alpha(\frac{1}{p}-\frac{1}{2}),\alpha}_{p}(\mathbb{R})$, and hence one has  $\mathcal{H}^p_{m_{s,\alpha}}(\cG_\rho,\tilde \cG_\rho) = M^{ s + \alpha(\frac{1}{p}-\frac{1}{2}),\alpha}_{p}(\mathbb{R})$. In particular, by Theorem \ref{banachfr}, $\cG_\rho$ is an atomic decomposition and a Banach frame for $M^{ s + \alpha(\frac{1}{p}-\frac{1}{2}),\alpha}_{p}(\mathbb{R})$.
Recall here that $g^\rho = g - g_\rho$.
Since  $\|\L\langle f, \mathcal{P}_j^\alpha M_{p_\alpha(j)} D_{s_\alpha^{-1}(j)} T_{a k} \tilde g_\rho\rangle \R_{j,k \in \mathbb{Z}} \|_{\ell^p_{m_{s,\alpha}}(\mathbb{Z}^2)}\leq B \|f\|_{M^{ s + \alpha(\frac{1}{p}-\frac{1}{2}),\alpha}_{p}}$, where $B>0$ is uniform with respect to $\rho$, to show (b) it is sufficient to verify that for all $\varepsilon>0$ there exists $\rho_0>0$ such that for all $\rho \geq \rho_0$, $a=a(\rho)\asymp \rho^{-1}$ and for any finite sequence $c=(c_{j,k})_{j,k \in \mathbb{Z}}$ of scalars
\begin{equation}
\label{f34}
\| \sum_{j,k \in \mathbb{Z}} c_{j,k}  M_{p_\alpha(j)} D_{s_\alpha^{-1}(j)} T_{a k} g^\rho\|_{M^{ s + \alpha(\frac{1}{p}-\frac{1}{2}),\alpha}_{p}} \leq \varepsilon \|c\|_{\ell^p_{m_{s,\alpha}}}. 
\end{equation}
Then one can apply \cite[Theorem 2.3]{CH}. By the equivalence of norms \eqref{f33} for some fixed $\rho^* \geq 1$ one has
\begin{eqnarray*}
&&\left \| \sum_{j,k \in \mathbb{Z}} c_{j,k}  M_{p_\alpha(j)} D_{s_\alpha^{-1}(j)} T_{a k} g^\rho \right \|_{M^{ s + \alpha(\frac{1}{p}-\frac{1}{2}),\alpha}_{p}}\\ &\asymp & \left \| \left (\sum_{j,k \in \mathbb{Z}} c_{j,k}  \langle M_{p_\alpha(j)} D_{s_\alpha^{-1}(j)} T_{a k} g^\rho, \mathcal{P}_i^\alpha M_{p_\alpha(i)} D_{s_\alpha^{-1}(i)} T_{a h} \tilde g_{\rho^*} \rangle \right)_{i,h}\right \|_{\ell^p_{m_{s,\alpha}}}.
\end{eqnarray*}
By an application of Lemma \ref{lemma32} (c) and Proposition \ref{boundedness} one has 
\begin{eqnarray*}
&&\| \sum_{j,k \in \mathbb{Z}} c_{j,k}  M_{p_\alpha(j)} D_{s_\alpha^{-1}(j)} T_{a k} g^\rho\|_{M^{ s + \alpha(\frac{1}{p}-\frac{1}{2}),\alpha}_{p}}\\ &\lesssim &D_\rho \left \| \left (\sum_{j,k \in \mathbb{Z}}  \left (1 + |j-i| \right)^{\frac{1}{2}\left(\frac{\alpha}{(1-\alpha)}\gamma_t-\gamma_f' \right)}   \left( 1+  \max\{s_\alpha(i),s_\alpha(j)\} |k s_\alpha(j)^{-1} - h s_\alpha(i)^{-1}| \right)^{-\frac{\gamma_t}{2}}|c_{j,k}|  \right)_{i,h}\right\|_{\ell^p_{m_{s,\alpha}}}\\
&\lesssim& D_\rho \|c\|_{\ell^p_{m_{s,\alpha}}}.
\end{eqnarray*}
Since $D_\rho \rightarrow 0$ for $\rho \rightarrow +\infty$, one shows \eqref{f34}.\\

Let us show  (c). First we have to observe that, by a direct computation, the operator $S$ defined by
\begin{eqnarray}
S(c) = \sum_{j,k \in\mathbb{Z}} c_{j,k}  \mathcal{P}_j^\alpha M_{p_\alpha(j)} D_{s_\alpha^{-1}(j)} T_{a k} \tilde g_\rho
\end{eqnarray}
is bounded from $\ell^p_{m_{s,\alpha}}$ into $M^{ s + \alpha(\frac{1}{p}-\frac{1}{2}),\alpha}_{p}(\mathbb{R})$ uniformly with respect to $\rho>0$ and $a=a(\rho)$.
\begin{eqnarray*}
\|S(c)\|_{M^{ s + \alpha(\frac{1}{p}-\frac{1}{2}),\alpha}_{p}(\mathbb{R})}^p &=& \sum_{j' \in \mathbb{Z}} \left \|  \mathcal{P}_{j'}^\alpha (\sum_{j,k \in\mathbb{Z}} c_{j,k}  \mathcal{P}_j^\alpha M_{p_\alpha(j)} D_{s_\alpha^{-1}(j)} T_{a k} \tilde g_\rho )\right \|_p^p (1+|p_\alpha(j'))^{sp + \alpha(1-\frac{p}{2})}\\
&=&    \sum_{j' \in \mathbb{Z}} \left \|\sum_{j \in \mathbb{Z}^d} \mathcal{P}_{j'}^\alpha  \mathcal{P}_j^\alpha M_{p_\alpha(j)} D_{s_\alpha^{-1}(j)} \left(\sum_{k \in\mathbb{Z}} c_{j,k}   T_{a k} \tilde g_\rho \right) \right \|_p^p (1+|p_\alpha(j'))^{sp + \alpha(1-\frac{p}{2})}\\
&\lesssim&  \sum_{j' \in \mathbb{Z}}  \sum_{j: P_j^\alpha P_{j'}^\alpha \neq 0} \left \|  \mathcal{P}_{j'}^\alpha  \mathcal{P}_j^\alpha M_{p_\alpha(j)} D_{s_\alpha^{-1}(j)} \left( \sum_{k \in\mathbb{Z}}c_{j,k}  T_{a k} \tilde g_\rho \right) \right \|_p^p (1+|p_\alpha(j'))^{sp + \alpha(1-\frac{p}{2})}\\
&\lesssim& \sum_{j' \in \mathbb{Z}}  \sum_{j: P_j^\alpha P_{j'}^\alpha \neq 0} s_\alpha(j)^{\frac{p-2}{2}}\left \|  \sum_{k \in\mathbb{Z}}c_{j,k}  T_{a k} \tilde g_\rho \right \|_p^p (1+|p_\alpha(j'))^{sp + \alpha(1-\frac{p}{2})}
\end{eqnarray*}
\begin{eqnarray*}
&=& \sum_{j' \in \mathbb{Z}}  \sum_{j: P_j^\alpha P_{j'}^\alpha \neq 0} s_\alpha(j)^{\frac{p-2}{2}}\left \|  a \sum_{k \in\mathbb{Z}}c_{j,k}  T_{a k} \tilde g \right \|_p^p (1+|p_\alpha(j'))^{sp + \alpha(1-\frac{p}{2})}\\
&\lesssim& \sum_{j' \in \mathbb{Z}}  \sum_{j: P_j^\alpha P_{j'}^\alpha \neq 0} s_\alpha(j)^{\frac{p-2}{2}}  a \|(c_{j,k})_k\|_{\ell^p(\mathbb{Z})}^p (1+|p_\alpha(j'))^{sp + \alpha(1-\frac{p}{2})}\\
&\lesssim&   a \sum_{j' \in \mathbb{Z}}  \|(c_{j',k})_k\|_{\ell^p(\mathbb{Z})}^p m_{s,\alpha}(j')^p  \leq \|c\|_{\ell_{m_{s,\alpha}}^p}^p.
\end{eqnarray*}
The first and the last inequality holds because the sum over $\{j: P_j^\alpha P_{j'}^\alpha \neq 0\}$ is uniformly finite.
Moreover, we have used $\|\sum_k d_k T_{ak} \tilde g\|_p \lesssim a^{1/p-1} \|d\|_{\ell^p}$.
Then it is sufficient to observe as before that for all 
$\varepsilon>0$ there exists $\rho_0>0$ such that for all $\rho\geq\rho_0$, and $a=a(\rho)\asymp \rho^{-1}$ and  for all $f \in M^{ s + \alpha(\frac{1}{p}-\frac{1}{2}),\alpha}_{p}$ 
\begin{equation}
\label{f35}
\| (\langle f, M_{p_\alpha(j)} D_{s_\alpha^{-1}(j)} T_{a k} g^\rho \rangle)_{j,k} \|_{ \ell^p_{m_{s,\alpha}}} \leq \varepsilon \|f\|_{ M^{ s + \alpha(\frac{1}{p}-\frac{1}{2}),\alpha}_{p}}, 
\end{equation}
since $f= \sum_{i,h} c_{j,k}  M_{p_\alpha(i)} D_{s_\alpha^{-1}(i)} T_{a h} g_{\rho^*} $ with $\|c\|_{\ell^p_{m_{s,\alpha}}} \asymp \|f\|_{M^{ s + \alpha(\frac{1}{p}-\frac{1}{2}),\alpha}_{p}}$.
These conditions are then enough to apply \cite[Theorem 2.2]{CH}. 
\end{proof}


\begin{rems} 1. In the proof of the previous theorem we have used the rather general and abstract results  \cite[Theorem 2.2, Theorem 2.3]{CH}. Of course, it is possible to keep the argument more concrete. In particular, it is not difficult to show that the operator
$$
S_\rho f = \sum_{j,k} \langle f, \mathcal{P}_j^\alpha M_{p_\alpha(j)} D_{s_\alpha^{-1}(j)} T_{a k} \tilde g_\rho \rangle M_{p_\alpha(j)} D_{s_\alpha^{-1}(j)} T_{a k} g,
$$
is bounded on $M^{ s + \alpha(\frac{1}{p}-\frac{1}{2}),\alpha}_{p}$. By Lemma \ref{lemma32} and Proposition \ref{boundedness} one can even  show that  for $\rho>0$ large enough
$$
\| I-S_\rho\|_{M^{ s + \alpha(\frac{1}{p}-\frac{1}{2}),\alpha}_{p} \rightarrow M^{ s + \alpha(\frac{1}{p}-\frac{1}{2}),\alpha}_{p}} <1.
$$
This implies that $S_\rho$ is boundedly invertible for $\rho>0$ large enough and that for all $f \in M^{ s + \alpha(\frac{1}{p}-\frac{1}{2}),\alpha}_{p}$ one has the unconditional convergent expansion
$$
f = S_\rho S_\rho^{-1} f = \sum_{j,k} \langle f, (S_\rho^{-1})^* \mathcal{P}_j^\alpha M_{p_\alpha(j)} D_{s_\alpha^{-1}(j)} T_{a k} \tilde g_\rho \rangle M_{p_\alpha(j)} D_{s_\alpha^{-1}(j)} T_{a k} g.
$$

2. The assumption $\mathcal{F} g \neq 0$ on $\Omega_0=[-1,1]$ is technical and it is essentially a non-vanishing condition. We expect that it can be removed.

3. Theorem \ref{chamod} is a generalization of \cite[Theorem 13.5.3]{Gr} and \cite[Theorem 5.2]{gro1} (see also \cite{FoG}), corresponding to the case $\alpha=0$, where Gabor frame characterizations of modulation spaces have been given. We conjecture that Theorem \ref{chamod} can be formulated for the case $\alpha \rightarrow 1$ to characterize inhomogeneous Besov spaces $B_p^{s-1/p-1/2}(\mathbb{R})$. Since $\lim_{\alpha \rightarrow 1} m_{s,\alpha}(j,k) = e^{s |j|}$, we expect that the extension of our theory to the case $\alpha \rightarrow 1$ should involve \emph{exponentially localized} frames as described in the previous Remark, see also \cite{gro1}. Interesting results in this direction have been suggested by Cordero and Gr\"ochenig in \cite{cg} for the wavelet frame characterization of \emph{homogeneous} Besov spaces.  


4. Theorem \ref{chamod} extends to the frame characterization of $M_{p,q}^{s,\alpha}(\mathbb{R})$ for $p \neq q$, just considering $\ell^{p,q}_{m_{s,\alpha}}$ spaces instead of $\ell^p_{m_{s,\alpha}}$. In fact, similarly to Proposition \ref{boundedness} and by applying standard arguments of complex interpolation of mixed norm sequence spaces \cite{BL}, matrices in $\cA_{\alpha,(1-\alpha)\gamma,\eta}$ are also bounded on $\ell^{p,q}_{m_{s,\alpha}}(\mathbb{Z}^2)$ for suitable $s$. 

5. Lemma \ref{lemma32} is strongly dependent on the particular geometry of the $\alpha$-covering determined by $(p_\alpha,s_\alpha)$ on the real line. We expect that the approach illustrated in this paper can be useful also for a frame characterization of $M_{p,q}^{s,\alpha}(\mathbb{R}^d)$ for $d>1$, with major technical difficulties.
\end{rems}

\subsection{$\alpha$-modulation spaces and time-frequency transforms}

In several relevant contributions, for example \cite{SAG,BI,CF,FF,Folland,Fo2,HL,HL2,HN,T1,T2}, an ``intermediate'' time-frequency transform between wavelet and short time Fourier transform is considered. 

Assume $\alpha \in [0,1]$ and $c>0$. For any $g \in L^2(\mathbb{R}) \backslash \{0\}$ and for $f \in L^2(\mathbb{R})$
 we define the \emph{flexible Gabor-wavelet transform} (or {\it $\alpha$-transform}) by
\begin{eqnarray}
        V_g^\alpha(f)(x,\omega) &:=& \langle f, T_x M_\omega D_{c(1+|\omega|)^{-\alpha}} g \rangle\\
        &=& \int_{\mathbb{R}} f(t) \overline{T_x M_\omega D_{c(1+|\omega|)^{-\alpha}} g(t)} dt, \quad x,\omega \in \mathbb{R}.
\end{eqnarray}
The transform can naturally extend to distributions whenever $g \in \mathcal{S}(\mathbb{R})$. 
For $\alpha=0$ the transform $V_g^\alpha$ coincides with the well-known short time Fourier transform, while for $\alpha=1$ it is a slight modification of the wavelet transform.
In particular, the intermediate case $\alpha=1/2$ is the Fourier-Bros-Iagolnitzer transform \cite{BI}.
In \cite[Theorem 4.4]{HN} Holschneider and Nazaret proved a characterization of $L^2$-Sobolev spaces by pull back techniques based on $\alpha$-transforms. For a suitable choice of $g \in \mathcal{S} \backslash \{0\}$ (for example the Gaussian) one has 
\begin{equation}
f \in H^s(\mathbb{R}) \text{ if and only if } V_g^\alpha (f) \in L^2_{m}(\mathbb{R}^2),
\end{equation}
where $m(x,\omega)=(1+|\omega|)^{s}$, $x, \omega \in \mathbb{R}$. In particular the following equivalence of norms holds
\begin{equation}
\|f\|_{H^s(\mathbb{R})} \asymp \|V_g^\alpha (f) \|_{L^2_m(\mathbb{R}^2)}, \text{ for all } f \in H^s(\mathbb{R}).
\end{equation}
Inspired by this characterization, they introduce a more general class of Banach spaces \cite[Definition 4.7]{HN}.
For a suitable choice of a Banach function space $B$ on the time-frequency plane $\mathbb{R}^2$ one can define the space of distributions on $\mathbb{R}$ given by
\begin{equation}
\mathbb{B}(\mathbb{R}):=\{f \in\mathcal{S}'(\mathbb{R}): \quad V_g^\alpha(f) \in B\}, 
\end{equation} 
endowed with the retract norm
\begin{equation}
\|f\|_{\mathbb{B}(\mathbb{R})} = \|V_g^\alpha(f)\|_B.
\end{equation}

A similar approach can be found in \cite[Section 4.6]{HL2} where generalizations of modulation spaces are introduced by Hogan and Lakey.

We want to observe here that, for the choice of $B$ as a certain weighted Lebesgue mixed norm $L^{p,q}$ space, the corresponding $\mathbb{B}(\mathbb{R})$ space is an $\alpha$-modulation space.
In fact, since $f \in M_{p,q}^{s,\alpha}(\mathbb{R})$ if and only if $\mathcal{F} f \in D(\cI_\alpha,\mathcal{F} L^p,\ell^q_{w_s})$, the \emph{decomposition space} subordinate to the covering $\cI_\alpha$, with local component $\mathcal{F} L^p$, and global component $\ell^q_{w_s}(\cI_\alpha)$ (see \cite{F,FG,G} for details), by an application of \cite[Theorem 4.3]{F} one can show the following 
\begin{tm}
Assume $s \in \mathbb{R}$, $\alpha \in [0,1]$, and $1 \leq p,q <  \infty$. For a suitable band-limited $g \in \mathcal{S}(\mathbb{R}) \backslash \{0\}$
\begin{equation}
M_{p,q}^{s+\alpha(1/q-1/2),\alpha}(\mathbb{R})=\{f \in \mathcal{S}'(\mathbb{R}): V_g^\alpha(f) \in L^{p,q}_{m}(\mathbb{R}^2)\}.
\end{equation}  
Moreover the norm of $M_{p,q}^{s+\alpha(1/q-1/2),\alpha}(\mathbb{R})$ can be equivalently expressed by
\begin{equation}
\|f\|_{M_{p,q}^{s+\alpha(1/q-1/2),\alpha}(\mathbb{R})} \asymp \L\int_{\mathbb{R}} \L \int_{\mathbb{R}} |V_g^\alpha(f)(x,\omega)|^p dx \R^{q/p} (1+|\omega|)^{s q} d\omega \R^{1/q},
\end{equation}  
for all $f \in M_{p,q}^{s+\alpha(1/q-1/2),\alpha}(\mathbb{R})$.
For $p \cdot q=\infty$ the usual modifications apply.
\end{tm}
A detailed discussion on the relations between continuous and discrete characterization of $\alpha$-modulation spaces will be given elsewhere in the context of recent generalizations of the coorbit space theory \cite{DST,DST2,FR}.

\subsection{Equivalence of frames and $\alpha$-modulation spaces}
As we have seen, qualities of frames can be observed by studying their associated Banach spaces.
Therefore, the ``differences'' between associated Banach spaces can be considered a ``measure'' of the different analysis that two frames perform.
The results in this paper can be interpreted as a qualitative study of the ``degree of difference'' of the analysis performed by Gabor and wavelet frames (Fig. 1). 

Let us conclude recalling in the following some of the relevant results related to inclusions of $\alpha$-modulations spaces, investigated by Gr\"obner \cite{G}:
\begin{tm}
If $1 \leq p,q \leq \infty$, $s \in \mathbb{R}$ and $0 \leq \alpha_1 < \alpha_2 \leq 1$ then
\begin{equation}
        M_{p,q}^{s',\alpha_2}(\mathbb{R}) \subset M_{p,q}^{s,\alpha_1}(\mathbb{R}), \quad s' = s + \frac{( \alpha_2 - \alpha_1)}{q}
\end{equation}
\begin{equation}
M_{p,q}^{s,\alpha_1}(\mathbb{R}) \subset M_{p,q}^{s',\alpha_2}(\mathbb{R}), \quad s' = s -  ( 1- 1/q) (\alpha_2 - \alpha_1).
\end{equation}
In particular, for $\alpha_2=1$ and $\alpha_1=0$,
\begin{equation}
B_{p,q}^{s+1/q}(\mathbb{R}) \subset M_{p,q}^{s}(\mathbb{R}).
\end{equation}
\end{tm}

\end{document}

%% file: macros.tex
\newtheorem{tm}{Theorem}[section]
\newtheorem{lemma}[tm]{Lemma}
\newtheorem{prop}[tm]{Proposition}

 \theoremstyle{definition}
 \newtheorem{definition}{Definition}
\newtheorem{exmp}{Example}
\newtheorem{exmps}{Examples}

\newcommand{\beqa}{\begin{eqnarray*}}
\newcommand{\eeqa}{\end{eqnarray*}}

\newcommand{\field}[1]{\mathbb{#1}}
\newcommand{\bR}{\field{R}}        
\newcommand{\bZ}{\field{Z}}        
 \def\cF{\mathcal{F}}              
 \def\cS{\mathcal{S}}
 
 \def\cH{\mathcal{H}}
 \def\cB{\mathcal{B}}
 
 \def\cG{\mathcal{G}}

 \def\cA{\mathcal{A}}
 
 \def\cI{\mathcal{I}}

 \def\cX{\mathcal{X}}
 
\def\cP{\mathcal{P}}

\def\rd{\bR^d}

\def\zd{\bZ^d}

\def\L{\left(}
\def\R{\right)}

\def\<{\left<}
\def\>{\right>}

\def\mv1{M_v^1}


%% file: framesAMcorrected2.bbl
\begin{thebibliography}{9}
\bibitem{SAG} {S. T. Ali, J. P. Antoine, J. P. Gazeau}, \textit{Coherent States, Wavelets and their Generalizations}, Springer-Verlag, 2000.
\bibitem{cg} E. Cordero, K. Gr\"ochenig, Localization of frames II, Appl. Comp. Harmon. Anal., {\bf 17}, no. 1, 2004, pag. 29-47.
\bibitem{bchl} R. Balan, P. Casazza, C. Heil, Z. Landau, Density, redundancy, and localization of frames, preprint, 2003.
\bibitem{BL} J. Bergh, J. L\"ofstr\"om, {\it Interpolation spaces. An introduction}, Springer-Verlag, Berlin 1976. Grundlehren der Mathematischen Wissenschaften No. 223.
\bibitem{B} L. Borup, Pseudodifferential operators on $\alpha$-modulation spaces,  J. Func. Spaces and Appl., {\bf 2}, no. 2, May 2004.
\bibitem{BI} J. Bros, D. Iagolnitzer, Support essentiel et structure analytique des distributions, in Seminaire Goulaouic-Lions-Schwartz, exp. no. 18, 1975.
\bibitem{chr} {O. Christensen}, {\it An Introduction to Frames and Riesz Bases}, Birkh\"auser, 2003.
\bibitem{CH} O. Christensen, C. Heil, Perturbation of Banach frames and atomic decompositions, Math. Nachr. {\bf 185}, 1997, pag. 33-47. 
\bibitem{CF}  A. Cordoba, C. Fefferman, Wave packets and Fourier integral operators, Comm. Partial Diff. Eq.,{\bf 3}, 1978, pag. 979-1005.
\bibitem{DFR} S. Dahlke, M. Fornasier, and T. Raasch,  Adaptive frame methods for elliptic operator equations, Bericht 2004-3, Fachbereich Mathematik und Informatik, Philipps-Universit\"at Marburg, Germany, 2004.
\bibitem{DST} {S. Dahlke, G. Steidl, and G. Teschke}, {Coorbit spaces and Banach frames on homogeneous spaces with applications to analyzing functions on spheres}, Adv. Comp. Math., {\bf 21}, no. 1, 2004, 147-180.
\bibitem{DST2} {S. Dahlke, G. Steidl, and G. Teschke}, {Weighted coorbit spaces and Banach frames on homogeneous spaces}, to appear in J. Four. Anal. Appl.
\bibitem{D1}  {I. Daubechies},  {Wavelets, time-frequency localization and signal analysis}, IEEE Trans. Inf. Th., {\bf 36}, 1990, pag. 961-1005.
\bibitem{D2}  {I. Daubechies}, \textit{Ten Lectures on Wavelets}, SIAM, 1992.
\bibitem{DGM} I. Daubechies, A. Grossmann, Y. Meyer, Painless nonorthogonal expansions, J. Math Phys. {\bf 27}, no. 5,  1986, pag. 1271-1283.
\bibitem{DT}  {L. Daudet, B. Torresani}, {Hybrid representations for audiophonic signal encoding}, Sign. Proc. {\bf 82}, no. 11, 2002, pag. 1585-1617.
\bibitem{DS}  {R. J. Duffin, A. C. Schaeffer}, A class of nonharmonic Fourier series, Trans. Amer. Math. Soc., {\bf 72}, 1952, pag. 341-366.
\bibitem{F}  {H. G. Feichtinger}, {Banach spaces of distributions defined by decomposition methods II}, Math. Nachr., {\bf 132}, 1987, pag. 207-237.
\bibitem{F4} {H. G. Feichtinger}, {Atomic characterization of modulation spaces through Gabor-type representations}, Proc. Conf. Constr. Function Theory, Rocky Mountain J. Math. {\bf 19}, 1989, pag. 113-126.
\bibitem{FF}  H. G. Feichtinger, M. Fornasier, Flexible Gabor-wavelets atomic decompositions for $L^2$-Sobolev spaces, to appear in Annali di Matematica Pura e Applicata.
\bibitem{FG} {H. G. Feichtinger, P. Gr\"obner}, {Banach spaces of distributions defined by decomposition methods I},  Math. Nachr., {\bf 123}, 1985, pag. 97-120.
\bibitem{FG1} {H. G. Feichtinger, K. Gr\"ochenig}, {A unified approach to atomic decomposition via integrable group representations}, Springer Lect. Notes Math., {\bf  1302}, 1988.
\bibitem{FG2}  {H. G. Feichtinger, K. Gr\"ochenig}, {Banach spaces related to integrable group representations and their atomic decompositions I}, J. Funct. Anal., {\bf  86}, 1989, pag. 307-340.
\bibitem{FG3} {H. G. Feichtinger, K. Gr\"ochenig}, {Banach spaces related to integrable group representations and their atomic decompositions II}, Monatsh. f. Math. {\bf 108}, 1989, pag. 129-148.
\bibitem{FG5} {H. G. Feichtinger, K. Gr\"ochenig}, {Irregular sampling theorems and series expansions of band-limited functions}, J. Math. Anal. Appl., {\bf 167}, 1992, pag. 530-556.
\bibitem{FS} {H. G. Feichtinger, T. Strohmer (Eds.)}, \textit{Gabor Analysis and Algorithms}, Birkh\"auser, 1998.
\bibitem{FS1} {H. G. Feichtinger, T. Strohmer (Eds.)}, {\it Advances in Gabor Analysis}, Birkh\"auser, 2003.
\bibitem{Folland} {G. B. Folland}, \textit{Harmonic Analysis in Phase Space}, Princeton Univ. Press, no. 122 in Annals Math. Studies, 1989.
\bibitem{Fo} {M. Fornasier}, {Decompositions of Hilbert spaces: local construction of global frames}, Proc. of ``Constructive Theory of Functions 2002'', Varna 2002, (B. Bojanov, Ed.), DARBA, Sofia, 2003,  pag. 255-281.
\bibitem{Fo1} {M. Fornasier}, {Quasi-orthogonal decompositions of frames}, J. Math. Anal. Appl., {\bf 289}, no. 1 , 2004, pag. 180-199.
\bibitem{Fo2} {M. Fornasier}, \textit{Constructive Methods for Numerical Applications in Signal Processing and Homogenization Problems}, Ph.D. thesis, University of Padova, 2002.
\bibitem{FoG} M. Fornasier, K. Gr\"ochenig, Intrinsic localization of frames, Preprint 8/2004, Dipartimento di Metodi e Modelli Matematici per le Scienze Applicate, Universit\`a di Roma ``La Sapienza'', 2004.
\bibitem{FR} M. Fornasier, H. Rauhut, Continuous frames, function spaces, and the discretization problem, preprint, 2004.
\bibitem{FJ} {M. Frazier, B. Jawerth}, {Decomposition of Besov spaces}, {\it Indiana Univ. Math. J.} {\bf 34}, 1985,  pag. 777-799.
\bibitem{G} {P. Gr\"obner}, \textit{Banachr\"aume glatter Funktionen and Zerlegungmethoden}, Ph.D. thesis, University of Vienna, 1992.
\bibitem{gro} {K. Gr\"ochenig}, {Describing functions: atomic decompositions versus frames}, Monatsh. Math. {\bf 112}, 1991,  pag. 1-41.
\bibitem{Gr} {K. Gr\"ochenig}, \textit{Foundation of Time-Frequency Analysis}, Birkh\"auser Verlag, 2001.
\bibitem{gro1} {K. Gr\"ochenig}, {Localization of frames, Banach frames, and the invertibility of the frame operator}, J. Four. Anal. Appl. {\bf 10}, no. 2, 2004, pag. 105-132. 
\bibitem{GL} K.~Gr{\"o}chenig, M.~Leinert, Symmetry of matrix algebras and symbolic calculus for infinite matrices, preprint 2003.
\bibitem{HLW} {E. Hernandez, D. Labate, G. Weiss},{ A unified characterization of reproducing systems generated by a finite family II}, J. Geom. Anal. {\bf 12}, no. 4, 2002, pag. 615-662.  
\bibitem{HLWW} {E. Hernandez, D. Labate, G. Weiss, E. Wilson}, {Oversampling quasi affine frames and wave packets}, Appl. Comp. Harmon. Anal., {\bf 16}, no. 2, 2002, 111-147.
\bibitem{HL} {J.A. Hogan, J.D. Lakey}, {Extensions of the Heisenberg group by dilations and frames}, Appl. Comp. Harmon. Anal. {\bf 2}, 1995,  pag. 174-199.
\bibitem{HL2} {J.A. Hogan, J.D. Lakey}, {Embeddings and uncertainty principles for generalized modulation spaces}, Chapater 4 in ``Modern Sampling Theory: Mathematics and Applications'', J. J. Benedetto and P. J. S. G. Ferreira (Eds.), Birkh\"auser, Boston, 2000, pag. 73-105.
\bibitem{HN} {M. Holschneider, B. Nazareth}, {An interpolation family between Gabor and wavelet transformations. Application to differential calculus and construction of anisotropic Banach spaces}, Adv. In Partial Diff. Eq., ``Nonlinear Hyperbolic Equations, Spectral Theory, and Wavelets Transformations'' (Albeverio, Demuth, Schrohe, Schulze Eds.), Wiley 2003, pag. 363-394.
\bibitem{j} S. Jaffard, Propri\'et\'es des matrices ``bien localis\'ees'' pre\`e de leur diagonale et qualques applications. Ann. Inst. H. Poicar\'e Anal. Non Lin\'eaire, {\bf 7}, no. 5, 1990,   pag. 461-476.
 \bibitem{L}{D. Labate},{ A unified characterization of reproducing systems generated by a finite family}, J. Geom. Anal. {\bf 12}, no. 3, 2002,  pag. 469-491.
\bibitem{Le} P. G. Lemari\'e, Bases d'ondelettes sur les groupes de Lie stratifi\'e, Bulletin de la Societ\'e Math\'ematique de France, {\bf 177}, no. 2, 1989, pag. 213-232.
\bibitem{M} Y. Meyer, {\it Ondelettes et operateurs I, II, III}, Hermann, Paris, 1990-1991.
\bibitem{NB} M. Nielsen, L. Borup, Nonlinear approximation in $\alpha$-modulation spaces, preprint 2003.
\bibitem{PS} L. P\"aiv\"arinta, E. Somersalo, A generalization of the Calderon-Vaillancourt theorem to $L^p$ and $h^p$, Math. Nachr. {\bf 138}, 1988, pag. 145-156.
\bibitem{St} R.~Stevenson, Adaptive solution of operator equations using wavelet frames, SIAM J. Numer. Anal. {\bf 41}, no. 3, pag. 1074-1100.
\bibitem{T1} {B. Torresani}, {Wavelets associated with representations of the affine Weyl-Heisenberg group}, J. Math. Phys. {\bf 32}, 1991, pag. 1273-1279.
\bibitem{T2} {B. Torresani}, {Time-frequency representation: wavelet packets and optimal decomposition}, Ann. Inst. H. Poincar\'e {\bf 56}, 1992, pag. 215-234.
\bibitem{Tr} {H. Triebel}, \textit{Theory of Function Spaces}, Birkh\"auser, 1983.
\bibitem{Tr2} {H. Triebel}, \textit{Theory of Function Spaces II}, Birkh\"auser, 1992.





\end{thebibliography}
